\documentclass[review]{elsarticle}

\usepackage{amssymb}
\usepackage{graphicx}
\usepackage{graphics}

\usepackage{lineno,hyperref}
\modulolinenumbers[5]

\journal{Annales de l'Institut Henri Poincar\'e (C) Analyse Non
Lin\'eaire}

\begin{document}

\begin{frontmatter}

\title{Optimal homotopy perturbation method for nonlinear differential equations governing MHD Jeffery-Hamel flow with heat transfer problem}

\author[VMarinca]{Vasile Marinca}
\ead{vmarinca@mec.upt.ro}

\author[RDEne]{Remus-Daniel Ene\corref{mycorrespondingauthor}}
\cortext[mycorrespondingauthor]{Corresponding author.}
\ead{remus.ene@upt.ro}


\address[VMarinca]{University Politehnica Timi\c soara, Department of Mechanics and Vibration, Timi\c soara, 300222, Romania\\
Department of Electromechanics and Vibration, Center for Advanced
and Fundamental Technical Research, Romania Academy, Timi\c soara,
300223, Romania}
\address[RDEne]{University Politehnica Timi\c soara, Department of
Mathematics, Timi\c soara, 300006, Romania}



%
%

\begin{abstract}
In this paper, Optimal Homotopy Perturbation Method (OHPM) is
employed to determine an analytic approximate solutions for
nonlinear MHD Jeffery-Hamel flow and heat transfer problem. The
Navier-Stokes equations, taking into account Maxwell's
electromagnetism and heat transfer lead to two nonlinear ordinary
differential equations. The obtained results by means of OHPM show
a very good agreement in comparison with the numerical results and
with Homotopy Perturbation Method (HPM).
\end{abstract}

\begin{keyword}
optimal homotopy perturbation method \sep Jeffery-Hamel \sep
nonlinear ordinary differential equations.
\end{keyword}

\end{frontmatter}

\linenumbers

\section{Introduction}
\label{1}

\noindent \par The incompressible fluid flow with heat transfer is
one of the most applicable cases in various fields of engineering
due to it industrial applications. The problem of a viscous fluid
between two nonparallel walls meeting at a vertex and with a
source of sink at the vertex was pioneered by Jeffery
\cite{Jeffery1915}, and Hamel \cite{Hamel1916}. Later, the
Jeffery-Hamel problem have been studied by several researchers and
discussed in many textbooks and articles. A stationary problem
with a finite number of "outlets" to infinity in the form of
infinite sectors is considered by Rivkind and Solonnikov
\cite{Rivkind2000}. The problem of steady viscous flow in a
convergent channel is analyzed analytically and numerically for
small, moderately large and asymptotically large Reynolds numbers
over entire range of allowed convergence angles by Akulenko et al.
\cite{Akulenko2004}. The MHD Jeffery-Hamel problem is solved by
Makinde and Mhone \cite{Makinde2006} using a special type of
Hermite-Pad\'e approximation semi-numerical approach and by
Esmaili et al. \cite{Esmaili2008} by applying Adomian
decomposition method. The classical Jeffery-Hamel flow problem is
solved by Ganji et al. \cite{Ganji2009} by means of the
variational iteration method and homotopy perturbation method, and
by Joneidi et al. \cite{Joneidi2010}, by differential
transformation method, Homotopy Perturbation Method and Homotopy
Analysis Method. The classical Jeffery-Hamel problem was extended
in \cite{Moghimi2011} to include the effects of external magnetic
field in conducted fluid. Optimal homotopy asymptotic method is
applied by Marinca and Heri\c sanu \cite{Marinca12011} and by
Esmaeilpour and Ganji \cite{Esmaeilpour2010}. The effect of
magnetic field and nanoparticle on the Jeffery-Hamel flow are
studied in \cite{Sheikholeshami2012} and \cite{Rostami2014}.
Numerical treatment using stochastic algorithms is used by Raja
and Samar \cite{Raja2014}.

In general, the problems as Jeffery-Hamel flows and other fluid
mechanics problems are inherently nonlinear. Excepting a limited
number of these problems most do not have analytical solution. The
aim of this paper is to propose an accurate approach to the MHD
Jeffery-Hamel flow with heat transfer problem using an analytical
technique, namely OHPM \cite{Marinca22011}, \cite{Marinca32011},
\cite{Marinca42012}. Our approach does not require a small or
large parameter in the governing equations, is based on the
construction and determination of some auxiliary functions
combined with a convenient way to optimally control the
convergence of the solution.

\section{Problem statement and governing equations}
\label{2}

\noindent \par We consider a system of cylindrical coordinates
with a steady flow of an incompressible conducting viscous fluid
from a source or sink at channel walls lying in planes, with angle
$2 \alpha$, taking into account the effect of electromagnetic
induction, as shown in Fig. 1, and the heat transfer.\\


The
continuity equation, the Navier-Stokes equations and energy
equation in cylindrical coordinates, can be written as
\cite{Schlichting1979}, \cite{Ali2011}, \cite{Choi1979},
\cite{Turkylmazoglu2014}:\\
\begin{eqnarray}\label{JefferyHamel1}  
\frac{1}{r} \frac{\partial }{\partial r} \left(r u_r \right)+
\frac{1}{r} \frac{\partial }{\partial r} \left(r u_{\varphi}
\right) = 0
\end{eqnarray}
\begin{eqnarray}\label{JefferyHamel2}  
\begin{array}{ll}
u_r \frac{\partial u_r}{\partial r} + \frac{u_{\varphi}}{r}
\frac{\partial u_r}{\partial \varphi} - \frac{u_{\varphi}^2}{r} =
- \frac{1}{\rho} \frac{\partial  P}{\partial
r} 
+ \nu \left[ \frac{1}{r} \frac{\partial (r
\varepsilon_{rr})}{\partial r} + \frac{1}{r} \frac{\partial
\varepsilon_{r\varphi}}{\partial r} -
\frac{\varepsilon_{r\varphi}}{r} - \frac{\sigma B_0^2 }{\rho r^2}
u_{r} \right]
\end{array}
\end{eqnarray}
\begin{eqnarray}\label{JefferyHamel3}  
\begin{array}{ll}
u_r \frac{\partial u_{\varphi}}{\partial r} +
\frac{u_{\varphi}}{r} \frac{\partial u_{\varphi}}{\partial
\varphi} - \frac{u_{\varphi} u_r}{r} = - \frac{1}{\rho r}
\frac{\partial P}{\partial
{\varphi}} 
+ \nu \left[ \frac{1}{r^2} \frac{\partial (r \varepsilon_{r
{\varphi}})}{\partial r} + \frac{1}{r} \frac{\partial
\varepsilon_{\varphi \varphi}}{\partial {\varphi}} -
\frac{\varepsilon_{r \varphi}}{r} - \frac{\sigma B_0^2 }{\rho r^2}
u_{\varphi} \right]
\end{array}
\end{eqnarray}
\begin{eqnarray}\label{JefferyHamel4}  
\begin{array}{ll}
u_r \frac{\partial T}{\partial r} = \frac{k}{\rho c_p}
\left(\nabla^2 T \right) + \frac{\nu}{c_p} \left[ 2 \left(
\left(\frac{\partial u_r}{\partial r} \right)^2 + \left(\frac{
u_r}{ r} \right)^2 \right) + \left( \frac{1}{ r} \frac{\partial
u_r}{\partial r} \right)^2 \right] + \frac{\sigma B_0^2 }{\rho
r^2} u_{r}^2
\end{array}
\end{eqnarray}
where $\rho$ is the fluid density, $P$ is the pressure, $\nu$ is
the kinematic viscosity, $T$ is the temperature, $k$ is the
thermal conductivity, $c_p$ is the specific heat at constant
pressure, $\sigma$ is the electrical conductivity, $B_0$ is the
induced magnetic field and the stress components are defined as\\
\begin{eqnarray}\label{JefferyHamel5}  
\varepsilon_{rr}=2\frac{\partial u_r}{\partial r } - \frac{2}{3}
div \bar{u}
\end{eqnarray}
\begin{eqnarray}\label{JefferyHamel6}  
\varepsilon_{\varphi \varphi} = \frac{2}{r} \frac{\partial
u_{\varphi}}{\partial \varphi } + \frac{2 u_r}{r} - \frac{2}{3}
div \bar{u}
\end{eqnarray}
\begin{eqnarray}\label{JefferyHamel7}  
\varepsilon_{r \varphi} = \frac{2}{r} \frac{\partial u_r}{\partial
\varphi } + 2 \frac{\partial}{\partial r}
\left(\frac{u_{\varphi}}{r}\right)
\end{eqnarray}

By considering the velocity field is only along radial direction
i.e. $u_{\varphi}=0$ and substituting Eqs.
(\ref{JefferyHamel5})-(\ref{JefferyHamel7}) into Eqs.
(\ref{JefferyHamel2}) and (\ref{JefferyHamel3}), the continuity,
Navier-Stokes and energy equations become:\\
\begin{eqnarray}\label{JefferyHamel8}  
\frac{1}{r} \frac{\partial }{\partial r} \left(r u_r \right) =  0
\end{eqnarray}
\begin{eqnarray}\label{JefferyHamel9}  
\begin{array}{ll}
u_r \frac{\partial u_r}{\partial r} = - \frac{1}{\rho}
\frac{\partial  P}{\partial r} 
+ \nu \left( \nabla^2 u_r - \frac{u_r}{r^2} - \frac{\sigma B_0^2
}{\rho r^2} u_{r} \right)
\end{array}
\end{eqnarray}
\begin{eqnarray}\label{JefferyHamel10}  
\begin{array}{ll}
- \frac{1}{\rho r} \frac{\partial P}{\partial \varphi} 
+ \frac{2 \nu}{r^2} \frac{\partial u_r }{\partial \varphi} = 0
\end{array}
\end{eqnarray}

The relevant boundary conditions, due to the symmetry assumption
at the channel centerline are as follows:
\begin{eqnarray}\label{JefferyHamel11}  
\begin{array}{ll}
\frac{\partial u_r }{\partial \varphi} = \frac{\partial T
}{\partial \varphi} = 0, \ \ \ u_r=\frac{u_c}{r} \ \ \ \textrm{at}
\ \ \varphi=0
\end{array}
\end{eqnarray}
and at the plates making body of channel:
\begin{eqnarray}\label{JefferyHamel12}  
\begin{array}{ll}
u_r=0, \ \ \  T = \frac{T_c}{r^2} \ \ \ \textrm{at} \ \ \varphi=
\alpha
\end{array}
\end{eqnarray}
where $u_c$ and $T_c$ are the centerline rate of movement and the
constant wall temperature, respectively.

From the continuity equation (\ref{JefferyHamel8}), one can get\\
\begin{eqnarray}\label{JefferyHamel13}  
\begin{array}{ll}
r u_r= f(\varphi)
\end{array}
\end{eqnarray}
where $f(\varphi)$ is an arbitrary function of $\varphi$ only.

By integrating Eq. (\ref{JefferyHamel10}) it holds that\\
\begin{eqnarray}\label{JefferyHamel14}  
\begin{array}{ll}
P(r, \varphi) = \frac{2 \rho \nu}{r^2} f(\varphi) + \rho g(r)
\end{array}
\end{eqnarray}
in which $g(r)$ is an arbitrary function of $r$ only.

Now, defining the dimensionless parameters:\\
\begin{eqnarray}\label{JefferyHamel15}  
\begin{array}{ll}
\eta = \frac{\varphi}{\alpha}, \ \ \ F(\eta) = \frac{f(
\varphi)}{u_c}, \ \ \  \theta(\eta) = r^2 \frac{T}{T_c}
\end{array}
\end{eqnarray}
where $T_c$ is the ambient temperature, and substituting these
into Eqs. (\ref{JefferyHamel4}) and (\ref{JefferyHamel9}) and then
eliminating the pressure term, one can put:\\
\begin{eqnarray}\label{JefferyHamel16}  
\begin{array}{ll}
F'''+2\alpha Re F F' + (4-H)\alpha^2 F' = 0
\end{array}
\end{eqnarray}
\begin{eqnarray}\label{JefferyHamel17}  
\begin{array}{ll}
\theta'' + 2 \alpha \left(2\alpha + Re Pr F \right) \theta + \beta
Pr \left[\left(H+4\alpha^2\right) F^2 + {F'}^2 \right] = 0
\end{array}
\end{eqnarray}
subject to the boundary conditions\\
\begin{eqnarray}\label{JefferyHamel18}  
\begin{array}{ll}
F(0)=1, \ \ \ F'(0)=0, \ \ \ F(1)=0
\end{array}
\end{eqnarray}
\begin{eqnarray}\label{JefferyHamel19}  
\begin{array}{ll}
\theta(1)=0, \ \ \ \theta'(0)=0
\end{array}
\end{eqnarray}
where $Re=\frac{\alpha u_c}{\nu}$ is the Reynolds number,
$H=\sqrt{\frac{\sigma B_0^2}{\rho \nu}}$ is the Hartmann number,
$Pr=\frac{\nu c_p}{k \rho}$, $\beta=\frac{u_c}{c_p}$ and prime
denotes derivative with respect to $\eta$.

\section{Basic ideas of optimal homotopy perturbation method}
\label{3}

\noindent \par To explain the ideas of the optimal homotopy
perturbation method, consider the non-linear differential equation
\begin{eqnarray}\label{JefferyHamel20}  
L\left[u,u',u'',u''', \eta \right] + g(\eta) +
N\left[u,u',u'',u''', \eta \right] =0 
\end{eqnarray}
that is subject to the initial / boundary condition:
\begin{eqnarray}\label{JefferyHamel21}  
B\left(u, \frac{\partial u}{\partial \eta}\right)=0, \ \ \ \eta
\in \Gamma
\end{eqnarray}
where $L$ is a linear operator, $g$ is a known function, $N$ a
nonlinear operator, $B$ is a boundary operator and $\Gamma$ is the
boundary of the domain of interest \cite{Marinca22011},
\cite{Marinca32011}, \cite{Marinca42012}. We construct the
homotopy \cite{He2000}\\
\begin{eqnarray}\label{JefferyHamel22}  
{\mathcal{H}}\left(u, p \right) = L\left(u,u',u'',u''', \eta
\right) + g(\eta) + p  N\left( u,u',u'',u''', \eta  \right)=0
\end{eqnarray}
for Eq. (\ref{JefferyHamel20}), where $p$ is the homotopy
parameter, $p \in \left[0, \ 1 \right]$. From Eq.
(\ref{JefferyHamel22}) one gets:\\
\begin{eqnarray}\label{JefferyHamel23}  
\begin{array}{ll}
{\mathcal{H}}\left(u, 0 \right) = L\left(u,u',u'',u''', \eta
\right) + g(\eta) =0\\
{\mathcal{H}}\left(u, 1 \right) = L\left(u,u',u'',u''', \eta
\right) + g(\eta) + N\left( u,u',u'',u''', \eta  \right)=0
\end{array}
\end{eqnarray}

Assuming that the approximate analytical solution of the
second-order can be expressed in the form\\
\begin{eqnarray}\label{JefferyHamel24}  
\bar{u}(\eta) = u_0+p u_1+p^2 u_2
\end{eqnarray}
and expanding the nonlinear operator $N$ in series, with respect
to the parameter $p$, we have:\\
\begin{eqnarray}\label{JefferyHamel25}  
\begin{array}{ll}
N\left( {\bar{u}},{\bar{u}}',{\bar{u}}'',{\bar{u}}''', \eta
\right) = N\left( u_0, u'_0,u''_0,u'''_0, \eta  \right) + p \left[
u_1 N_{\bar{u}}\left( u_0, u'_0,u''_0,u'''_0, \eta  \right) + \right. \\
+ u'_1 N_{{\bar{u}}'}\left( u_0, u'_0,u''_0,u'''_0, \eta \right) +
u''_1 N_{{\bar{u}}''}\left( u_0, u'_0,u''_0,u'''_0, \eta \right) + \\
\left. + u'''_1 N_{{\bar{u}}'''}\left( u_0, u'_0,u''_0,u'''_0,
\eta \right) \right] + p^2 \left( u_2 N_{\bar{u}} + u'_2
N_{{\bar{u}}'} + ...  \right)
\end{array}
\end{eqnarray}
where $F_{\bar{u}} = \frac{\partial F}{\partial \bar{u}}$. By
introducing a number of unknown auxiliary functions $H_i(\eta,
C_k)$, $i=0, \ 1, \ 2, \ ...$ that depend on the variable $\eta$
and some parameters $C_k$, $k=  1, \ 2, \ ..., \ s$, we can construct a new homotopy:\\
\begin{eqnarray}\label{JefferyHamel26}  
\begin{array}{ll}
{\mathcal{H}}\left(\bar{u}, p \right) = L\left(
{\bar{u}},{\bar{u}}',{\bar{u}}'',{\bar{u}}''', \eta \right) +
g(\eta)+  p
H_0(\eta, C_k) N\left( u_0, u'_0,u''_0,u'''_0, \eta  \right) + \\
p^2 \left[ H_1(\eta, C_k) u_1 N_{\bar{u}}\left( u_0 \right) +
H_2(\eta, C_k) u'_1 N_{{\bar{u}}'}\left( u_0 \right) + H_3(\eta,
C_k) u''_1 N_{{\bar{u}}''}\left( u_0 \right) + \right. \\
\left. + H_4(\eta, C_k) u'''_1 N_{{\bar{u}}'''}\left( u_0 \right)
\right] + p^2 \left[ H_5(\eta, C_k) u_2 N_{\bar{u}}\left( u_0
\right) + H_6(\eta, C_k) u'_2 N_{{\bar{u}}'}\left( u_0 \right) +
... \right]
\end{array}
\end{eqnarray}

Equating the coefficients of like powers of $p$, yields the linear
equations:\\
\begin{eqnarray}\label{JefferyHamel27}  
L\left[u_0, u'_0,u''_0,u'''_0, \eta  \right] + g(\eta)=0, \quad
B\left(u_0, \frac{\partial u_0 }{\partial \eta}\right)=0
\end{eqnarray}
\begin{eqnarray}\label{JefferyHamel28}  
\begin{array}{ll}
L\left(u_1 \right) + H_0(\eta, C_k) N\left( u_0,
u'_0,u''_0,u'''_0, \eta  \right) =0, \quad B\left(u_1,
\frac{\partial u_1 }{\partial \eta}\right)=0
\end{array}
\end{eqnarray}
\begin{eqnarray}\label{JefferyHamel29}  
\begin{array}{ll}
L\left(u_2 \right) + H_1(\eta, C_k) u_1 N_{\bar{u}}\left( u_0,
u'_0,u''_0,u'''_0, \eta \right) + H_2(\eta, C_k) u'_1
N_{{\bar{u}}'}\left( u_0 \right) + \\
+ H_3(\eta, C_k) u''_1 N_{{\bar{u}}''}\left( u_0 \right) +
H_4(\eta, C_k) u'''_1 N_{{\bar{u}}'''}\left( u_0 \right) = 0,
\quad B\left(u_2, \frac{\partial u_2 }{\partial \eta}\right)=0
\end{array}
\end{eqnarray}

The functions $H_i(\eta, C_k)$, $i=0, \ 1, \ 2, \ ...$ are not
unique and can be chosen such that the products $H_i \cdot u_j
N_u$ and $u_j N_u$ are of the same form. In this way, a maximum of
only two iterations are required to achieve accurate solutions.

The unknown parameters $C_k$, $k= 1, \ 2, \ ..., \ s$ which appear
in the functions $H_i(\eta, C_k)$ can be determined optimally by
means of the least-square method, collocation method, the weighted
residuals, the Galerkin method, and so on.

In this way the solution of Eq. (\ref{JefferyHamel20}) subject to
the initial / boundary condition (\ref{JefferyHamel21}) can be
readily determined. It follows that the basic ideas of our
procedure are the construction of a new homotopy
(\ref{JefferyHamel26}), the auxiliary functions $H_i$ with
 parameters $C_k$ that can be determined optimally leading to the
 conclusion that the convergence of the approximate solutions can
 be easily controlled.

\section{Application of OHPM to the MHD Jeffery-Hamel flow and heat transfer problem}
\label{4}

\noindent \par Let us present the approximate analytic expressions
of $f(\eta)$ and $\theta(\eta)$ from Eqs.
(\ref{JefferyHamel16})-(\ref{JefferyHamel19}) by means of OHPM.

For Eqs. (\ref{JefferyHamel16}) and (\ref{JefferyHamel18}), the
linear operator is chosen as $L\left(F\right) = F'''$, while the
nonlinear operator is defined as $N\left(F\right) = 2\alpha F F' +
(4-H)\alpha^2 F'$, $g(\eta)=0$. The initial approximation $F_0$ is
obtained from equation (\ref{JefferyHamel27})\\
\begin{eqnarray}\label{JefferyHamel30}  
F_0'''=0, \ \ \ \quad F_0(0)=1, \ \ F'_0(0)=0, \ \ F_0(1)=0
\end{eqnarray}

The solution of Eq. (\ref{JefferyHamel30}) is hence
\begin{eqnarray}\label{JefferyHamel31}  
\begin{array}{ll}
F_0(\eta) = 1-\eta^2.
\end{array}
\end{eqnarray}

On the other hand, from Eq. (\ref{JefferyHamel16}), one obtains
\begin{eqnarray}\label{JefferyHamel32}  
N_{F}\left(F\right)=2\alpha Re F', \ \ \ \quad
N_{F'}\left(F\right)=2\alpha Re F + (4-H)\alpha^2.
\end{eqnarray}

By substituting Eq. (\ref{JefferyHamel31}) into the nonlinear
operator $N$ and into Eq. (\ref{JefferyHamel32}) one retrieves:\\
\begin{eqnarray}\label{JefferyHamel33}  
\begin{array}{ll}
N\left(F_0\right)= 2A\eta^2-2(A+B)\eta, \ \ \
N_{F}\left(F_0\right)=-2A \eta, \\
N_{F'}\left(F_0\right)= -A
\eta^2 +A+B
\end{array}
\end{eqnarray}
where $A=2\alpha Re$, $B=(4-H)\alpha^2$.

Eq. (\ref{JefferyHamel28}) becomes
\begin{eqnarray}\label{JefferyHamel34}  
F'''_1 + H_0(\eta, C_k)\left[ 2A\eta^2-2(A+B)\eta \right]=0, \quad
F_1(0)=F'_1(0)=F_1(1)=0
\end{eqnarray}

We choose $H_0(\eta, C_k)=-60 C_1$ where $C_1$ is an unknown
parameter and from Eq. (\ref{JefferyHamel34}) we obtain\\
\begin{eqnarray}\label{JefferyHamel35}  
\begin{array}{ll}
F_1(\eta) = 2A C_1 \eta^5 -5 (A+B)C_1 \eta^4 + (3A+5B)C_1 \eta^2
\end{array}
\end{eqnarray}

Eq. (\ref{JefferyHamel29}) can be written in the form
\begin{eqnarray}\label{JefferyHamel36}  
\begin{array}{ll}
F'''_2 + H_1(\eta, C_k) (-2A\eta) F_1 + H_2(\eta, C_k) (-A\eta^2
+A+B) F'_1 =0\\
F_2(0)=F'_2(0)=F_2(1)=0
\end{array}
\end{eqnarray}

In this case we choose
$$H_1(\eta, C_k)= \frac{1}{2 A} \left(C_2 \eta^2 + C_3 \eta + C_4 + \frac{C_5}{\eta} \right), \quad H_2(\eta, C_k)= \frac{C_6}{2} + \frac{C_7}{\eta}$$
such that the solution of Eq. (\ref{JefferyHamel36}) is given by\\
\begin{eqnarray}\label{JefferyHamel37}  
\begin{array}{ll}
F_2(\eta) = \frac{AC_1 C_2}{495} \eta^{11} + \frac{2 A C_1 C_3 -5
(A+B) C_1 C_2}{720} \eta^{10} + \\
+ \frac{2 A C_1 C_4 -5 (A+B) C_1 C_3 + 5 A^2 C_1 C_6}{504}
\eta^{9} + \\
+ \frac{(3A + 5B) C_1 C_2 - 5 (A+B) C_1 C_4 + 2 A C_1 C_5 - 10
(A^2+ A B) C_1 C_6}{336} \eta^{8} +\\
+ \frac{(3A + 5B) C_1 C_3 - 5 (A+B) C_1 C_5 - 5
(A^2+ A B) C_1 C_6 + 2A C_1 C_7}{210} \eta^{7} +\\
+ \frac{(3A + 5B) C_1 C_4 +
(13 A^2+ 25 A B + 10 B^2) C_1 C_6 - 5(A+B)C_1 C_7}{120} \eta^{6} +\\
+ \frac{(3A + 5B) C_1 C_5}{60} \eta^{5} - \frac{(3 A^2+ 8 A B + 5
B^2) C_1 C_6 + (3A+5B)C_1 C_7}{24} \eta^{4} + M \eta^{2}
\end{array}
\end{eqnarray}
where
\begin{eqnarray*} 
\begin{array}{ll}
M= -C_1 C_2 \left[ \left(\frac{1}{495} - \frac{1}{144} +
\frac{1}{112} \right) A + \left(\frac{5}{336} - \frac{1}{144}
\right) B \right]-\\
 -C_1 C_3 \left[ \left(\frac{1}{360} - \frac{5}{504} +
\frac{1}{70} \right) A + \left(\frac{1}{42} - \frac{5}{504}
\right) B \right]-\\
-C_1 C_4 \left[ \left(\frac{1}{252} - \frac{5}{336} + \frac{1}{40}
\right) A  + \left( \frac{1}{24} - \frac{5}{336} \right) B
\right]-C_1 C_5 \left(\frac{9 A}{280} + \frac{5
B}{84} \right) - \\
-C_1 C_6 \left[ \left(\frac{5}{504} - \frac{5}{168} -
\frac{5}{210} - \frac{1}{60} \right) A^2 -
\left(\frac{5}{168} + \frac{5}{210} + \frac{1}{8} \right) A B -  \frac{1}{8} B^2 \right] - C_1 C_7 \left( \frac{13 A}{140} + \frac{B}{6} \right)  
\end{array}
\end{eqnarray*}

For $p=1$ into Eq. (\ref{JefferyHamel24}), we obtain the
second-order approximate solution, using and Eqs.
(\ref{JefferyHamel31}), (\ref{JefferyHamel35}) and
(\ref{JefferyHamel37}):\\
\begin{eqnarray}\label{JefferyHamel38}  
\begin{array}{ll}
\bar{F}(\eta) = F_0(\eta) + F_1(\eta) + F_2(\eta) = \frac{AC_1
C_2}{495} \eta^{11} + \frac{2 A C_1 C_3 -5
(A+B) C_1 C_2}{720} \eta^{10} + \\
+ \frac{2 A C_1 C_4 -5 (A+B) C_1 C_3 + 5 A^2 C_1 C_6}{504}
\eta^{9} + \\
+ \frac{(3A + 5B) C_1 C_2 - 5 (A+B) C_1 C_4 + 2 A C_1 C_5 - 10
(A^2+ A B) C_1 C_6}{336} \eta^{8} +\\
+ \frac{(3A + 5B) C_1 C_3 - 5 (A+B) C_1 C_5 - 5
(A^2+ A B) C_1 C_6 + 2A C_1 C_7 }{210} \eta^{7} +\\
+ \frac{(3A + 5B) C_1 C_4 +
(13 A^2+ 25 A B + 10 B^2) C_1 C_6 - 5(A+B)C_1 C_7 }{120} \eta^{6} 
+ \left(2A C_1 + \frac{(3A + 5B) }{60} C_1 C_5 \right) \eta^{5} -
\\
- \left[ \frac{(3 A^2+ 8 A B + 5 B^2)C_1 C_6 + (3A+5B)C_1 C_7
}{24}  + (5 A + 5 B)
C_1 \right] \eta^{4} + \left[ \left( 3A+5B\right) + M-1 \right]\eta^{2} 
\end{array}
\end{eqnarray}

Now, we present the approximate analytic solution for Eqs.
(\ref{JefferyHamel17}) and (\ref{JefferyHamel19}). The linear and
nonlinear operators and the function $g$ are, respectively\\
\begin{eqnarray}\label{JefferyHamel39}  
\begin{array}{ll}
L(\theta) = \theta'', \ \ \ g(\eta)=-1, \\
N(\theta) = 1 + 4
\alpha^2 \theta + 2 \alpha Re Pr F \theta + \beta Pr \left[ (1+4
\alpha^2)F^2 + {F'}^2 \right]
\end{array}
\end{eqnarray}

The Eq. (\ref{JefferyHamel27}) becomes
\begin{eqnarray}\label{JefferyHamel40}  
\begin{array}{ll}
\theta''_0 -1=0, \ \ \ \quad \theta_0(1)=0, \ \ \ \theta'_0(0)=0
\end{array}
\end{eqnarray}

Eq. (\ref{JefferyHamel40}) has the solution
\begin{eqnarray}\label{JefferyHamel41}  
\begin{array}{ll}
\theta_0(\eta)=\frac{1}{2}(1-\eta^2).
\end{array}
\end{eqnarray}

From Eq. (\ref{JefferyHamel39}) it follows that
\begin{eqnarray}\label{JefferyHamel42}  
\begin{array}{ll}
N_{\theta}(\theta) = 4\alpha^2 + 2 \alpha Re Pr F, \ \ \ \quad
N_{{\theta}'}(\theta) =0.
\end{array}
\end{eqnarray}

By substituting Eq. (\ref{JefferyHamel41}) into Eqs.
(\ref{JefferyHamel39}) and (\ref{JefferyHamel42}), one gets
respectively\\
\begin{eqnarray}\label{JefferyHamel43}  
\begin{array}{ll}
N(\theta_0) = C-2D\eta^2+E\eta^4, \ \ \ N_{\theta}(\theta_0) =
L+K\eta^2
\end{array}
\end{eqnarray}
where
\begin{eqnarray}\label{JefferyHamel44}  
\begin{array}{ll}
C=1+2\alpha^2+\alpha Re Pr + 4 \beta \alpha^2 Pr + Pr H \beta\\
D=\alpha^2+2 \alpha Re Pr + 2 \beta Pr ( 2 \alpha^2 -1 ) + Pr H
\beta\\
E=\alpha Re Pr + 4 \beta \alpha^2 Pr + Pr H \beta\\
L=4\alpha^2 +2 \alpha Re Pr, \quad \quad \quad K= -2 \alpha Re Pr
\end{array}
\end{eqnarray}

Eq. (\ref{JefferyHamel28}) can put
\begin{eqnarray}\label{JefferyHamel45}  
\begin{array}{ll}
\theta''_1+h_0(\eta, C_8) (C-2D\eta^2+E\eta^4)=0, \ \ \ \quad
\theta_1(1)=\theta'_1(0)=0
\end{array}.
\end{eqnarray}

Choosing $h_0(\eta, C_8) = -30 C_8$ into Eq.
(\ref{JefferyHamel45}), one obtains\\
\begin{eqnarray}\label{JefferyHamel46}  
\begin{array}{ll}
\theta_1(\eta) = C_8 \left[ 15 C (\eta^2-1) -5D (\eta^4-1) +
E(\eta^6-1) \right]
\end{array}.
\end{eqnarray}

Eq. (\ref{JefferyHamel29}) can be written in the form
\begin{eqnarray}\label{JefferyHamel47}  
\begin{array}{ll}
\theta''_2+ h_1(\eta, C_k) (L+K\eta^2) \theta_1 = 0, \ \ \ \quad
\theta_2(1)=\theta'_2(0)=0
\end{array}
\end{eqnarray}
and therefore it is natural to choose the auxiliary function $h_1$
as\\
\begin{eqnarray*} 
\begin{array}{ll}
h_1(\eta, C_k) = \frac{1}{C_8} \left(C_9+C_{10} \eta +
C_{11}\eta^2 + C_{12} \eta^3 + C_{13} \eta^4\right)
\end{array}
\end{eqnarray*}

From Eq. (\ref{JefferyHamel47}), it can shown that
\begin{eqnarray}\label{JefferyHamel48}  
\begin{array}{ll}
\theta_2(\eta) = - L (15 C -5D +E) \left[ \frac{1}{2}C_9 (\eta^2
-1) + \frac{1}{6} C_{10} (\eta^3 -1)\right] + \\
+ \left[ 15 L C C_9 - (15 C -5D +E)(K C_9 + L C_{11}) \right]
\frac{\eta^4 -1}{12} + \\
+ \left[ 15 L C C_{10} - (15 C -5D +E)(K C_{10} + L C_{12})
\right]
\frac{\eta^5 -1}{20} + \\
+ \left[ (15C K -5 D K) C_9 + 15 L C C_{11} - \right.\\
\left. - (15 C -5D +E)(K
C_{11} + L C_{13}) \right] \frac{\eta^6 -1}{30} + \\
+ \left[ (15C K -5 D L) C_{10} + 15 L C C_{12} - K (15 C -5D +E) C_{12} \right] \frac{\eta^7 -1}{42} + \\
+ \left[ (L E -5 D K ) C_9 + (15 C K - 5 D L) C_{11} + 15 L C C_{13} - \right.\\
\left. - K (15 C -5 D +E) C_{13} \right] \frac{\eta^8 -1}{56} + \\
+ \left[ (L E -5 D K ) C_{10} + (15 C K - 5 D L) C_{12} \right] \frac{\eta^9 -1}{72} + \\
+ \left[ E K C_9 + (L E -5 D K ) C_{11} + (15 C K - 5 D L) C_{13} \right] \frac{\eta^{10} -1}{90} + \\
+ \left[ E K C_{11} + (L E -5 D K ) C_{13} \right] \frac{\eta^{12}
-1}{132} + E K C_{12} \frac{\eta^{13} -1}{156} + E K C_{13} \frac{\eta^{14} -1}{182} 
\end{array}
\end{eqnarray}

The second-order approximate solution of Eqs.
(\ref{JefferyHamel17}) and (\ref{JefferyHamel19}) is
\begin{eqnarray}\label{JefferyHamel49}  
\begin{array}{ll}
\bar{\theta}(\eta) = \theta_0(\eta) + \theta_1(\eta) +
\theta_2(\eta)
\end{array}
\end{eqnarray}

where $\theta_0$, $\theta_1$ and $\theta_2$ are given by Eqs.
(\ref{JefferyHamel41}), (\ref{JefferyHamel46}) and
(\ref{JefferyHamel48}) respectively.

\section{Numerical results}
\label{5}

\noindent \par In order to show the efficiency and accuracy of the
OHPM , we consider some cases for different values of the
parameters $\alpha$ and $H$. In all cases we consider $Re=50$,
$Pr=1$, $\beta = 3.492161428 \cdot 10^{-13}$.

\textbf{Case 5.1} Consider $\alpha = \frac{\pi}{24}$ and
 $H=0$. By means of the least-square method, the values of the parameters $C_i$, $i=1, \ 2, \ ..., \
 13$ are\\
$$C_1 = -0.000025737857, \ C_2 = 128165.4247848388, \ C_3 =
-360860.8730122449,$$ $$C_4 = 315974.73981422884, \ C_5 =
-20823.768289134576,$$ $$C_6 = -319.2089575339067, \ C_7 =
-62701.61063351235,$$ $$C_8 = -8.904447449602 \cdot 10^{-12}, \
C_9 = -2.9799239131772537 \cdot 10^{-12}, $$ $$C_{10} =
0.441808726864 \cdot 10^{-12}, \ C_{11}= -1.808778662628 \cdot
10^{-12}, $$ $$ C_{12} = 0.111699001609 \cdot 10^{-12}, \ C_{13} =
-5.293450955214 \cdot 10^{-12}$$

One get approximate solutions from Eqs. (\ref{JefferyHamel38}) and
(\ref{JefferyHamel49}) respectively:\\
\begin{eqnarray}\label{JefferyHamel50}  
\begin{array}{ll}
\bar{F}(\eta)= 1 - 2.3104494668 \eta^2 + 2.4868857696 \eta^4 +
 0.3531718162 \eta^5 - \\
 - 3.4153413394 \eta^6 +
 1.7515474936 \eta^7 + 1.2031805064 \eta^8 -\\
 -  1.6209066880 \eta^9 + 0.6391440832 \eta^{10} -
 0.0872321749 \eta^{11}
\end{array}
\end{eqnarray}
\begin{eqnarray}\label{JefferyHamel51}  
\begin{array}{ll}
\bar{\theta}(\eta)= \left[ -59.560673288998 (\eta^{2}-1) -
   9.116379402803 (\eta^{3}-1) - \right. \\
-142.753455187914 (\eta^{4}-1) +  44.843714524611 (\eta^{5}-1) +\\
+168.833421156648 (\eta^{6}-1) -  18.843109698135 (\eta^{7}-1) -\\
 -183.399486372212 (\eta^{8}-1) +   2.145923952284 (\eta^{9}-1) +\\
+ 122.114218185962 (\eta^{10}-1) - 36.681425691368 (\eta^{12}-1) -\\
\left. - 0.061343981363 (\eta^{13}-1) + 2.491809125293
(\eta^{14}-1) \right] \cdot 10^{-12}
\end{array}
\end{eqnarray}

\textbf{Case 5.2} For $\alpha = \frac{\pi}{24}$, $H=250$, the
parameters $C_i$ are:\\
$$C_{1} = -0.011565849071, \ C_{2} = 0.707159448177, \ C_{3} =
-290.324617452708,$$ $$C_{4} = 435.512207098156, \ C_{5} =
-98.780455208880,$$ $$C_{6} = -0.371954495218, \ C_{7} =
-82.314218258861,$$ $$C_{8} = -37.285501459040 \cdot 10^{-12}, \
C_{9} = -15.253404678160 \cdot 10^{-12},$$ $$C_{10} =
12.572809436948 \cdot 10^{-12}, \ C_{11} = -22.083271773195 \cdot
10^{-12},$$
$$C_{12} = 1.090339198400 \cdot 10^{-12}, \ C_{13} =
14.241169061094 \cdot 10^{-12}$$ 
and therefore the approximate solutions, (\ref{JefferyHamel38})
and (\ref{JefferyHamel49}) may be written
as:\\
\begin{eqnarray}\label{JefferyHamel52}  
\begin{array}{ll}
\bar{F}(\eta)= 1 - 1.638305627622 \eta^{2} + 1.206009669620
\eta^{4} +
 0.043648374556 \eta^{5} -\\
 - 1.078984595104 \eta^{6} +
 0.156296141929 \eta^{7} + 0.738925954400 \eta^{8} -\\
- 0.549972600570 \eta^{9} + 0.122598968736 \eta^{10} -
 0.000216285946 \eta^{11}
\end{array}
\end{eqnarray}
\begin{eqnarray}\label{JefferyHamel53}  
\begin{array}{ll}
\bar{\theta}(\eta)=  \left[  -1249.857282929460 (\eta^{2}-1) -
    9.116379418809 (\eta^{3}-1) -  \right.\\
 - 956.378989402938 (\eta^{4}-1) +1351.297888810558 (\eta^{5}-1)
  -\\
- 1018.834483075130 (\eta^{6}-1) - 646.214253910126 (\eta^{7}-1)
 +\\
+ 1232.314945834321 (\eta^{8}-1) + 129.244515099455 (\eta^{9}-1)
 -\\
-  589.041009734213 (\eta^{10}-1) +116.176788092689 (\eta^{12}-1)
  -\\
\left.  -  0.598803450274 (\eta^{13}-1) -  6.703807283538
(\eta^{14}-1) \right] \cdot 10^{-12}
\end{array}
\end{eqnarray}

\textbf{Case 5.3} For $\alpha = \frac{\pi}{24}$, $H=500$ we obtain:\\
\begin{eqnarray}\label{JefferyHamel54}  
\begin{array}{ll}
\bar{F}(\eta)= 1 - 1.1724778890 \eta^2 + 0.4686048815 \eta^4 -
 0.098961423266 \eta^5 - \\
 -0.3550702077 \eta^6 +
 0.788918793543 \eta^7 - 1.8412450049 \eta^8 +\\
+ 2.166714191416 \eta^9 - 1.2119558765 \eta^{10} +
 0.255472535060 \eta^{11}
\end{array}
\end{eqnarray}
\begin{eqnarray}\label{JefferyHamel55}  
\begin{array}{ll}
\bar{\theta}(\eta)=  \left[  -3025.399926380127 (\eta^{2}-1) -
    9.116379434815 (\eta^{3}-1) - \right.\\
- 3566.696711595512 (\eta^{4}-1) + 4417.102370290691 (\eta^{5}-1)
-\\
- 2124.393833887592 (\eta^{6}-1) - 2172.379321553169 (\eta^{7}-1)
+\\
+ 2808.629798261231 (\eta^{8}-1) +  471.928890384626 (\eta^{9}-1)
-\\
- 1267.595946519067 (\eta^{10}-1) + 212.551207395542 (\eta^{12}-1)
 -\\
\left. -    1.000321320275 (\eta^{13}-1) -  11.126037049812
(\eta^{14}-1) \right] \cdot 10^{-12}
\end{array}
\end{eqnarray}

\textbf{Case 5.4} For $\alpha = \frac{\pi}{24}$, $H=1000$, it holds that:\\
\begin{eqnarray}\label{JefferyHamel56}  
\begin{array}{ll}
\bar{F}(\eta)= 1 - 0.6223664982 \eta^2 - 0.1660109481 \eta^4 -
 0.2259232591 \eta^5 + \\
 +0.415889430872 \eta^6 -
 0.5758423027 \eta^7 + 0.198748794050 \eta^8 +\\
+ 0.0038449003 \eta^9 + 0.029311285675 \eta^{10} -
 0.0576514026 \eta^{11}
\end{array}
\end{eqnarray}
\begin{eqnarray}\label{JefferyHamel57}  
\begin{array}{ll}
\bar{\theta}(\eta)=  \left[  -8164.566371333924 (\eta^{2}-1) -
9.116379466826 (\eta^{3}-1) -  \right.\\
-12276.720695370957 (\eta^{4}-1) + 13537.207073598438 (\eta^{5}-1)
-\\
- 5217.969886112186 (\eta^{6}-1) - 6768.312414390118 (\eta^{7}-1)
 +\\
+ 8081.157802168993 (\eta^{8}-1) + 1537.402367705613 (\eta^{9}-1)
 -\\
- 3743.086450834657 (\eta^{10}-1) + 613.584536626787 (\eta^{12}-1)
 -\\
\left. -    1.303403201992 (\eta^{13}-1) -  31.505440255694
(\eta^{14}-1) \right] \cdot 10^{-12}
\end{array}
\end{eqnarray}

\textbf{Case 5.5} For $\alpha = \frac{\pi}{36}$, $H=0$, we obtain \\
\begin{eqnarray}\label{JefferyHamel58}  
\begin{array}{ll}
\bar{F}(\eta)= 1 - 1.7695466647 \eta^2 + 1.2754140854 \eta^4 +
 0.1103023597 \eta^5 - \\
 -1.1550256736 \eta^6 +
 0.4434051824 \eta^7 + 0.3358564399 \eta^8 -\\
- 0.3325192185 \eta^9 + 0.1045556849 \eta^{10} -
 0.0124421957 \eta^{11}
\end{array}
\end{eqnarray}
\begin{eqnarray}\label{JefferyHamel59}  
\begin{array}{ll}
\bar{\theta}(\eta)=  \left[  -111.208165621670 (\eta^{2}-1) -
6.895054111348 (\eta^{3}-1) - \right.\\
- 219.101666428999 (\eta^{4}-1) + 112.000412353392 (\eta^{5}-1) +\\
+ 120.895605108821 (\eta^{6}-1) -  47.637714830801 (\eta^{7}-1) -\\
- 104.355375516421 (\eta^{8}-1) +   6.984672282609 (\eta^{9}-1) +\\
+  75.737792150272 (\eta^{10}-1) - 24.982664253556 (\eta^{12}-1)
  -\\
\left. -   0.097447216002 (\eta^{13}-1) +  1.768576986049
(\eta^{14}-1) \right] \cdot 10^{-12}
\end{array}
\end{eqnarray}

\textbf{Case 5.6} If $\alpha = \frac{\pi}{36}$, $H=250$ then \\
\begin{eqnarray}\label{JefferyHamel60}  
\begin{array}{ll}
\bar{F}(\eta)= 1 - 1.5111014276 \eta^2 + 0.8612876996 \eta^4 +
 0.0128407190 \eta^5 - \\
 -0.5668759509 \eta^6 +
 0.0367598458 \eta^7 + 0.3198667899 \eta^8 -\\
- 0.1656253347 \eta^9 + 0.0019977968 \eta^{10} +
 0.0108498620 \eta^{11}
\end{array}
\end{eqnarray}
\begin{eqnarray}\label{JefferyHamel61}  
\begin{array}{ll}
\bar{\theta}(\eta)=  \left[  -5094.867741624686 (\eta^{2}-1) -
6.895054147361
(\eta^{3}-1) - \right.\\
- 4500.773385702808 (\eta^{4}-1) + 6249.12145418685 (\eta^{5}-1)
 -\\
- 4283.444052831423 (\eta^{6}-1) - 2834.31757828637 (\eta^{7}-1)
 +\\
+ 5102.615992457279 (\eta^{8}-1) + 528.701361603855 (\eta^{9}-1)
 -\\
- 2382.96513866767 (\eta^{10}-1) + 460.279590044557 (\eta^{12}-1)
 -\\
\left. -    2.850139405026 (\eta^{13}-1) - 26.598784857626
(\eta^{14}-1) \right] \cdot 10^{-12}
\end{array}
\end{eqnarray}

\textbf{Case 5.7} For $\alpha = \frac{\pi}{36}$, $H=500$ the approximate solutions are \\
\begin{eqnarray}\label{JefferyHamel62}  
\begin{array}{ll}
\bar{F}(\eta)= 1 - 1.2942214724 \eta^2 + 0.5334630869 \eta^4 +
 0.0019708504 \eta^5 - \\
 -0.3330602404 \eta^6 -
 0.0197377497 \eta^7 + 0.2108709627 \eta^8 -\\
- 0.1151093836 \eta^9 + 0.0147650020 \eta^{10} +
 0.0010589440 \eta^{11}
\end{array}
\end{eqnarray}
\begin{eqnarray}\label{JefferyHamel63}  
\begin{array}{ll}
\bar{\theta}(\eta)=  \left[  -10479.441564298402 (\eta^{2}-1) -
6.895054183374
(\eta^{3}-1) - \right.\\
- 9629.602717192312 (\eta^{4}-1) +12532.057010176612 (\eta^{5}-1)
 -\\
- 7708.77471180501 (\eta^{6}-1) -  5687.792565325205 (\eta^{7}-1)
 +\\
+ 9313.682156269786 (\eta^{8}-1) + 1063.274966246270 (\eta^{9}-1)
 -\\
- 4206.734533154069 (\eta^{10}-1) + 761.164805945711 (\eta^{12}-1)
 -\\
\left. -    5.659649234867 (\eta^{13}-1) -  42.382058179415
(\eta^{14}-1) \right] \cdot 10^{-12}
\end{array}
\end{eqnarray}

\textbf{Case 5.8} If $\alpha = \frac{\pi}{36}$, $H=1000$ then \\
\begin{eqnarray}\label{JefferyHamel64}  
\begin{array}{ll}
\bar{F}(\eta)= 1 - 0.9581900583 \eta^2 + 0.0913634769 \eta^4 +
 0.0030073285 \eta^5 - \\
 -0.1611067985 \eta^6 +
 0.0668486226 \eta^7 - 0.0696973759 \eta^8 +\\
+ 0.0454374262 \eta^9 - 0.0111066789 \eta^{10} -
 0.0065559425 \eta^{11}
\end{array}
\end{eqnarray}
\begin{eqnarray}\label{JefferyHamel65}  
\begin{array}{ll}
\bar{\theta}(\eta)=  \left[  -22457.262743604428 (\eta^{2}-1) -
6.895054255400
(\eta^{3}-1) - \right.\\
 22247.17539598977 (\eta^{4}-1) +  26605.166448464304 (\eta^{5}-1)
 -\\
- 14220.534699997283 (\eta^{6}-1) - 12131.493178637833
(\eta^{7}-1)
 +\\
+ 17921.197074275697 (\eta^{8}-1) +  2301.819882082563
(\eta^{9}-1)
 -\\
-  7836.564793794578 (\eta^{10}-1) + 1301.948445444968
(\eta^{12}-1)
  -\\
\left. -    11.186604132870 (\eta^{13}-1) -   68.576887703688
(\eta^{14}-1) \right] \cdot 10^{-12}
\end{array}
\end{eqnarray}

From the Tables 1-16 it is obvious that the second-order
approximate solutions obtained by OHPM are of a high accuracy in
comparison with homotopy perturbation method and with numerical
solution obtained by means of a fourth-order Runge-Kutta method in
combination with the shooting method using Wolfram Mathematica 6.0
software.

In Figs 2 and 3 are presented the effect of the Hartmann number on
the velocity profile for $Re=50$ and $\alpha=\frac{\pi}{24}$ and
$\alpha=\frac{\pi}{36}$ respectively. It is observe that velocity
increases with increasing of the Hartmann number for any value of
$\alpha$. The same effect of Hartmann number on the thermal
profile are presented in Figs 4 and 5 for $\alpha=\frac{\pi}{24}$
and $\alpha=\frac{\pi}{36}$ respectively. In this case, the
temperature decreases with increasing of the Hartmann number in
the both cases. The effect of the half angle $\alpha$ on the
velocity profile is presented in Figs 6-9. With an increasing
value of $\alpha$, velocity decreases for $H=0$ and $H=250$, but
increases for $H=500$ and $H=1000$. From the Figs 11-13, it is
interesting to remark that the temperature increases whereas the
half angle $\alpha$ increases. In all cases, the maximum of
temperature occurs near the walls for $H=0$ and precisely in wall
for $H \neq 0$, while the minimum occurs near the channel axis.\\

\section{Conclusions}
\label{6}

In present paper, the Optimal Homotopy Perturbation Method (OHPM)
is employed to propose a new analytic approximate solutions for
MHD Jeffery-Hamel flow with heat transfer problem. Our procedure
does not need restrictive hypotheses, is very rapid convergent
after only two iterations and the convergence of the solutions is
ensured in a rigorous way. The cornerstone of the validity and
flexibility of our procedure is the choice of the linear operator
and the optimal auxiliary functions which contribute to very
accurate solutions. The parameters which are involved in the
composition of the optimal auxiliary functions are optimally
identified via various methods in a rigorous way. Our technique is
very effective, explicit, easy to apply which proves that this
method is very efficient in practice.

\section*{Conflict of Interests}
\label{7}

The authors declare that there is no conflict of interests
regarding the publication of this paper.



\begin{thebibliography}{99}

\bibitem{Jeffery1915} G.~B.~ Jeffery, \emph{The two-dimensional steady motion of a viscous fluid}, 
Philosophical Magazine, 6(20), 1915, 455--465.

\bibitem{Hamel1916} G.~ Hamel, \emph{Spiralf$\ddot{o}$rmige bewegungen 
z$\ddot{a}$her flussigkeiten}, Jahresbericht der Deutschen
Mathematiker-Vereinigung, 25, 1916, 34--60.

\bibitem{Rivkind2000} L.~Rivkind, V.~A.~Solonnikov, \emph{Jeffery-Hamel asymptotics for 
steady state Navier-Stokes flow in domains with sector-like
outlets to infinity}, J. Math. Fluid Mech., 2, 2000, 324--352.

\bibitem{Akulenko2004} L.~D.~Akulenko, D.~V.~Georgevskii, S.~A.~Kumakshev, \emph{Solutions of the  
Jeffery-Hamel problem regularly extendable in the Reynolds
number}, Fluid Dynamics, 39(1), 2004, 12--28.

\bibitem{Makinde2006} O.~D.~Makinde, P.~Y.~Mhone, \emph{Hermite-Pad\'e approximation approach  
to MHD Jeffery-Hamel flows}, Applied Math. and Comput., 181, 2006,
966--972.

\bibitem{Esmaili2008} Q.~Esmaili, A.~Ramiar, E.~Alizadeh, D.~D.~Ganji, \emph{An approximation of the  
analytical solution of the Jeffery-Hamel flow by decomposition
method}, Phys. Lett., 372, 2008, 3434--3439.

\bibitem{Ganji2009} Z.~Z.~Ganji, D.~D.~Ganji, M.~Esmaeilpour, \emph{Study of 
nonlinear Jeffery-Hamel flow by He's semi-analytical methods and
comparison with numerical results}, Computers and Math with Appl.,
58, 2009, 2107--2116.

\bibitem{Joneidi2010} A.~A.~Joneidi, G.~Domairry, M.~Babaelahi, \emph{Three analytical    
methods applied to Jeffery-Hamel flow}, 
Commun. Nonlin. Sci. Numer. Simulat., 15, 2010, 3423--3434.

\bibitem{Moghimi2011} S.~M.~Moghimi, D.~D.~Ganji, H.~Bararnia, M.~Hosseini, M.~Jalbal, 
\emph{Homotopy perturbation method 
for nonlinear MHD Jeffery-Hamel problem}, Computers and Math. with
Appl., 61, 2011, 2213--2216.

\bibitem{Marinca12011} V.~Marinca, N.~Heri\c sanu, \emph{An optimal homotopy 
asymptotic approach applied to nonlinear MHD Jeffery-Hamel flow},
Mathematical Problems in Engineering, 2011, Article ID 169056, 16
pages.

\bibitem{Esmaeilpour2010} M.~Esmaeilpour, D.~D.~Ganji, \emph{Solution of the Jeffery-Hamel 
flow problem by optimal homotopy asymptotic method}, Comput. and
Math. Appl., 59, 2010, 3405--3411.

\bibitem{Sheikholeshami2012} M.~Sheikholeshami, D.~D.~Ganji, H.~R.~Ashorynejad, H.~B.~Rokni, 
\emph{Analytical investigation of Jeffery-Hamel flow with high
magnetic fiwld and nanoparticle by Adomian decomposition method},
Appl. Math. and Mech., Engl. ed., 33(1), 2012, 25--36.

\bibitem{Rostami2014} A.~K.~Rostami, M.~R.~Akbari, D.~D.~Ganji, S.~Heydari, \emph{Investigating 
Jeffery-Hamel flow with high-magnetic field and nanoparticles by
HPM and AGM}, Central Eur. J. of Eng., 4(4), 2014, 357--370.

\bibitem{Raja2014} M.~A.~Z.~Raja, R.~Samar, \emph{Numerical treatment of 
nonlinear MHD Jeffery-Hamel problems using stochastic algorithms},
Computers \& Fluids, 91, 2014, 111--115.

\bibitem{Marinca22011} V.~Marinca, N.~Heri\c sanu, \emph{Nonlinear Dynamical Systems in 
Engineering - Some Approximate Approaches}, Springer Verlag,
Heidelberg, 2011.

\bibitem{Marinca32011} V. Marinca, N. Heri\c sanu, \emph{Nonlinear dynamic analysis of an   
electrical machine rotor-bearing system by the optimal homotopy
perturbation method}, Computers and Math. with Appl., 61, 2011,
2019--2024.

\bibitem{Marinca42012} N. Heri\c sanu, V. Marinca, \emph{Optimal homotopy perturbation   
method for a non-conservative dynamical system of a rotating
electrical machine}, Z. Naturforsch, 67a, 2012, 509--516.

\bibitem{Schlichting1979} H.~Schlichting, \emph{Boundary-layer theory    
Mc Graw-Hill Book Company}, 1979.

\bibitem{Ali2011} F.~M.~Ali, R.~Nazar, N.~M.~Arifin, I.~Pop, \emph{MHD stagnation-point   
flow and heat transfer towards a stretching sheet with induced
magnetic field}, Appl. Math. Mech., 32, 2011, 409--418.

\bibitem{Choi1979} S.~H.~Choi, H.-E.~Wilhelm, \emph{Incompressible magnetohydrodynamic 
flow with heat transfer between inclined walls}, Physics of
Fluids, 22, 1979, 1073--1078.

\bibitem{Turkylmazoglu2014} M.~Turkylmazoglu, \emph{Extending the traditional 
Jeffery-Hamel flow to strechable convergent / divergent channels},
Computers \& Fluids, 100, 2014, 196--203.

\bibitem{He2000} J.~H.~He, \emph{A coupling  
method of homotopy technique and perturbation technique for
nonlinear problems}, Int. J. Non-Linear Mech., 35, 2000, 37--43.


\end{thebibliography}

%

\clearpage

\begin{table}[h!]\label{table1}
\textbf{Table 1.} Comparison between the HPM results
\cite{Moghimi2011}, OHPM results (\ref{JefferyHamel50}) and
numerical results for the velocity $F(\eta)$ for $\alpha =
\frac{\pi}{24}$ and $H=0$ \\
\begin{center}
\scriptsize{
\begin{tabular}{|c|c|c|c|c|}
  \hline
  $\eta$ & $F_{\textrm{HPM}}$ \cite{Moghimi2011} & $F_{\textrm{numeric}}$ & $ \bar{F}_{\textrm{OHPM}}$, Eq. (\ref{JefferyHamel50})& $\begin{array}{c}
                                                         \textrm{relative error}= \\
                                                          |F_{\textrm{numeric}}-\bar{F}_{\textrm{OHPM}}|
                                                          \end{array}$\\
                                                            \hline
0     & 1  &   1   &   1   &   0 \\  \hline 
0.1   & 0.9770711  &   0.9771426047   &   0.9771444959   &   1.8 $\cdot 10^{-6}$ \\  \hline 
0.2   & 0.9112020  &   0.9114792278   &   0.9114802054   &   9.7 $\cdot 10^{-7}$ \\  \hline 
0.3   & 0.8104115  &   0.8110052403   &   0.8110054657   &   2.2 $\cdot 10^{-7}$ \\  \hline 
0.4   & 0.6859230  &   0.6869148220   &   0.6869163034   &   1.4 $\cdot 10^{-6}$ \\  \hline 
0.5   & 0.5498427  &   0.5512883212   &   0.5512895302   &   1.2 $\cdot 10^{-6}$ \\  \hline 
0.6   & 0.4131698  &   0.4151088947   &   0.4151091740   &   2.7 $\cdot 10^{-7}$ \\  \hline 
0.7   & 0.2846024  &   0.2870546320   &   0.2870554998   &   8.6 $\cdot 10^{-7}$ \\  \hline 
0.8   & 0.1702791  &   0.1731221390   &   0.1731231521   &   1.01 $\cdot 10^{-6}$ \\  \hline 
0.9   & 0.0744232  &   0.0768715756   &   0.0768721058   &   5.3 $\cdot 10^{-7}$ \\  \hline 
1     & 0  &   0   &    0   &    0               \\  \hline 
    \end{tabular}}
\end{center}
\end{table}
\begin{table}[h!]\label{table2}
\textbf{Table 2.} Comparison between OHPM results
(\ref{JefferyHamel51}) and numerical results for the
temperature $\theta(\eta)$ for $\alpha = \frac{\pi}{24}$ and $H=0$ \\
\begin{center}
\scriptsize{
\begin{tabular}{|c|c|c|c|}
  \hline
  $\eta$ & ${\theta}_{\textrm{numeric}}$ & $ {\bar{{\theta}}}_{\textrm{OHPM}}$ from Eq. (\ref{JefferyHamel51})& $\begin{array}{c}
                                                         \textrm{relative error}= \\
                                                          |{{\theta}}_{\textrm{numeric}}-{\bar{{\theta}}}_{\textrm{OHPM}}|
                                                          \end{array}$\\
                                                            \hline
0     &   -9.134405103300 $\cdot 10^{-12}$   &   -9.134559900000 $\cdot 10^{-12}$   &   1.5 $\cdot 10^{-16}$ \\  \hline 
0.1   &   -8.553981690585 $\cdot 10^{-12}$   &   -8.561731325450 $\cdot 10^{-12}$   &   7.7 $\cdot 10^{-15}$ \\  \hline 
0.2   &   -7.029011445513 $\cdot 10^{-12}$   &   -7.029011445531 $\cdot 10^{-12}$   &   1.7 $\cdot 10^{-23}$ \\  \hline 
0.3   &   -4.968059593695 $\cdot 10^{-12}$   &   -4.959904647105 $\cdot 10^{-12}$   &   8.1 $\cdot 10^{-15}$ \\  \hline 
0.4   &   -2.830354806908 $\cdot 10^{-12}$   &   -2.830354806921 $\cdot 10^{-12}$   &   1.3 $\cdot 10^{-23}$ \\  \hline 
0.5   &   -1.003581673783 $\cdot 10^{-12}$   &   -1.015625108634 $\cdot 10^{-12}$   &   1.2 $\cdot 10^{-14}$ \\  \hline 
0.6    &   2.755715234490 $\cdot 10^{-13}$   &    2.755715234410 $\cdot 10^{-13}$   &   7.9 $\cdot 10^{-24}$ \\  \hline 
0.7    &   9.398183452650 $\cdot 10^{-13}$   &    9.682947296376 $\cdot 10^{-13}$   &   2.8 $\cdot 10^{-14}$ \\  \hline 
0.8    &   1.062337244632 $\cdot 10^{-12}$   &    1.062337244629 $\cdot 10^{-12}$   &   3.3 $\cdot 10^{-24}$ \\  \hline 
0.9    &   7.684282928672 $\cdot 10^{-13}$   &    6.524090716991 $\cdot 10^{-13}$   &   1.1 $\cdot 10^{-13}$ \\  \hline 
1      &   0   &    0                                 &   0         \\  \hline
    \end{tabular}}
\end{center}
\end{table}
\clearpage
\begin{table}[h!]\label{table3}
\textbf{Table 3.} Comparison between the HPM results
\cite{Moghimi2011}, OHPM results (\ref{JefferyHamel52}) and
numerical results for the velocity $F(\eta)$ for $\alpha =
\frac{\pi}{24}$ and $H=250$ \\
\begin{center}
\scriptsize{
\begin{tabular}{|c|c|c|c|c|}
  \hline
  $\eta$ &  $F_{\textrm{HPM}}$ \cite{Moghimi2011} &  $F_{\textrm{numeric}}$ & $ \bar{F}_{\textrm{OHPM}}$, Eq. (\ref{JefferyHamel52})& $\begin{array}{c}
                                                         \textrm{relative error}= \\
                                                          |F_{\textrm{numeric}}-\bar{F}_{\textrm{OHPM}}|
                                                          \end{array}$\\
                                                            \hline
  0   & 1  &   1   &   1   &   0 \\  \hline 
0.1   & 0.9837340  &   0.9837367791   &   0.9837369246   &   1.4 $\cdot 10^{-7}$ \\  \hline 
0.2   & 0.9363350  &   0.9363459948   &   0.9363459260   &   6.8 $\cdot 10^{-8}$ \\  \hline 
0.3   & 0.8616894  &   0.8617133581   &   0.8617132189   &   1.3 $\cdot 10^{-7}$ \\  \hline 
0.4   & 0.7653405  &   0.7653814753   &   0.7653813980   &   7.7 $\cdot 10^{-8}$ \\  \hline 
0.5   & 0.6533961  &   0.6534573087   &   0.6534570226   &   2.8 $\cdot 10^{-7}$ \\  \hline 
0.6   & 0.5314621  &   0.5315466609   &   0.5315462975   &   3.6 $\cdot 10^{-7}$ \\  \hline 
0.7   & 0.4038130  &   0.4039227354   &   0.4039224205   &   3.1 $\cdot 10^{-7}$ \\  \hline 
0.8   & 0.2728708  &   0.2729980317   &   0.2729975130   &   5.1 $\cdot 10^{-7}$ \\  \hline 
0.9   & 0.1389433  &   0.1390416079   &   0.1390411070   &   5.0 $\cdot 10^{-7}$ \\  \hline 
1     & 0  &   0   &    0   &    0 \\  \hline
    \end{tabular}}
\end{center}
\end{table}
\begin{table}[h!]\label{table4}
\textbf{Table 4.} Comparison between OHPM results
(\ref{JefferyHamel53}) and numerical results for the
temperature $\theta(\eta)$ for $\alpha = \frac{\pi}{24}$ and $H=250$ \\
\begin{center}
\scriptsize{
\begin{tabular}{|c|c|c|c|}
  \hline
  $\eta$ & ${\theta}_{\textrm{numeric}}$ & $ {\bar{{\theta}}}_{\textrm{OHPM}}$ from Eq. (\ref{JefferyHamel53})& $\begin{array}{c}
                                                         \textrm{relative error}= \\
                                                          |{{\theta}}_{\textrm{numeric}}-{\bar{{\theta}}}_{\textrm{OHPM}}|
                                                          \end{array}$\\
                                                            \hline
  0   &   -8.520036945014 $\cdot 10^{-10}$   &   -8.520036944914 $\cdot 10^{-10}$   &   9.9 $\cdot 10^{-21}$ \\  \hline 
0.1   &   -8.394493067234 $\cdot 10^{-10}$   &   -8.395974340239 $\cdot 10^{-10}$   &   1.4 $\cdot 10^{-13}$ \\  \hline 
0.2   &   -8.032504411694 $\cdot 10^{-10}$   &   -8.032504411616 $\cdot 10^{-10}$   &   7.8 $\cdot 10^{-21}$ \\  \hline 
0.3   &   -7.450795118612 $\cdot 10^{-10}$   &   -7.450297754118 $\cdot 10^{-10}$   &   4.9 $\cdot 10^{-14}$ \\  \hline 
0.4   &   -6.677062404515 $\cdot 10^{-10}$   &   -6.677062404494 $\cdot 10^{-10}$   &   2.1 $\cdot 10^{-21}$ \\  \hline 
0.5   &   -5.746573136109 $\cdot 10^{-10}$   &   -5.746735798562 $\cdot 10^{-10}$   &   1.6 $\cdot 10^{-14}$ \\  \hline 
0.6   &   -4.698328124476 $\cdot 10^{-10}$   &   -4.698328124489 $\cdot 10^{-10}$   &   1.3 $\cdot 10^{-21}$ \\  \hline 
0.7   &   -3.570882806332 $\cdot 10^{-10}$   &   -3.571240799348 $\cdot 10^{-10}$   &   3.5 $\cdot 10^{-14}$ \\  \hline 
0.8   &   -2.398013591946 $\cdot 10^{-10}$   &   -2.398013591937 $\cdot 10^{-10}$   &   9.8 $\cdot 10^{-22}$ \\  \hline 
0.9   &   -1.206024561793 $\cdot 10^{-10}$   &   -1.200670094842 $\cdot 10^{-10}$   &   5.3 $\cdot 10^{-13}$ \\  \hline 
 1    &   0   &    0                                 &   0 \\  \hline
    \end{tabular}}
\end{center}
\end{table}
\clearpage
\begin{table}[h!]\label{table5}
\textbf{Table 5.} Comparison between the HPM results
\cite{Moghimi2011}, OHPM results (\ref{JefferyHamel54}) and
numerical results for the velocity $F(\eta)$ for $\alpha =
\frac{\pi}{24}$ and $H=500$ \\
\begin{center}
\scriptsize{
\begin{tabular}{|c|c|c|c|c|}
  \hline
  $\eta$ & $F_{\textrm{HPM}}$ \cite{Moghimi2011} &  $F_{\textrm{numeric}}$ & $ \bar{F}_{\textrm{OHPM}}$, Eq. (\ref{JefferyHamel54})& $\begin{array}{c}
                                                         \textrm{relative error}= \\
                                                          |F_{\textrm{numeric}}-\bar{F}_{\textrm{OHPM}}|
                                                          \end{array}$\\
                                                            \hline
  0   & 1  &   1   &   1   &   0 \\  \hline 
0.1   & 0.9883197  &   0.9883196668   &   0.9883207994   &   1.1 $\cdot 10^{-6}$ \\  \hline 
0.2   & 0.9537955  &   0.9537952479   &   0.9538026351   &   7.3 $\cdot 10^{-6}$ \\  \hline 
0.3   & 0.8978515  &   0.8978510438   &   0.8978610430   &   9.9 $\cdot 10^{-6}$ \\  \hline 
0.4   & 0.8224699  &   0.8224690439   &   0.8224696000   &   5.5 $\cdot 10^{-7}$ \\  \hline 
0.5   & 0.7296817  &   0.7296803425   &   0.7296719365   &   8.4 $\cdot 10^{-6}$ \\  \hline 
0.6   & 0.6209748  &   0.6209728597   &   0.6209706881   &   2.1 $\cdot 10^{-6}$ \\  \hline 
0.7   & 0.4966644  &   0.4966617866   &   0.4966700156   &   8.2 $\cdot 10^{-6}$ \\  \hline 
0.8   & 0.3552115  &   0.3552079113   &   0.3552097261   &   1.8 $\cdot 10^{-6}$ \\  \hline 
0.9   & 0.1923821  &   0.1923775215   &   0.1923683642   &   9.1 $\cdot 10^{-6}$ \\  \hline 
1   &   0  &  0   &    0   &    0 \\  \hline
    \end{tabular}}
\end{center}
\end{table}
\begin{table}[h!]\label{table6}
\textbf{Table 6.} Comparison between OHPM results
(\ref{JefferyHamel55}) and numerical results for the
temperature $\theta(\eta)$ for $\alpha = \frac{\pi}{24}$ and $H=500$ \\
\begin{center}
\scriptsize{
\begin{tabular}{|c|c|c|c|}
  \hline
  $\eta$ & ${\theta}_{\textrm{numeric}}$ & $ {\bar{{\theta}}}_{\textrm{OHPM}}$ from Eq. (\ref{JefferyHamel55})& $\begin{array}{c}
                                                         \textrm{relative error}= \\
                                                          |{{\theta}}_{\textrm{numeric}}-{\bar{{\theta}}}_{\textrm{OHPM}}|
                                                          \end{array}$\\
                                                            \hline
  0   &   -1.783303641361 $\cdot 10^{-9}$   &   -1.783303641351 $\cdot 10^{-9}$   &   9.9 $\cdot 10^{-21}$ \\  \hline 
0.1   &   -1.753016328789 $\cdot 10^{-9}$   &   -1.753373570315 $\cdot 10^{-9}$   &   3.5 $\cdot 10^{-13}$ \\  \hline 
0.2   &   -1.666810282194 $\cdot 10^{-9}$   &   -1.666810282186 $\cdot 10^{-9}$   &   7.8 $\cdot 10^{-21}$ \\  \hline 
0.3   &   -1.531366088677 $\cdot 10^{-9}$   &   -1.531258066872 $\cdot 10^{-9}$   &   1.08 $\cdot 10^{-13}$ \\  \hline 
0.4   &   -1.356325163065 $\cdot 10^{-9}$   &   -1.356325163062 $\cdot 10^{-9}$   &   2.1 $\cdot 10^{-21}$ \\  \hline 
0.5   &   -1.152391353408 $\cdot 10^{-9}$   &   -1.152436518900 $\cdot 10^{-9}$   &   4.5 $\cdot 10^{-14}$ \\  \hline 
0.6   &   -9.301091200084 $\cdot 10^{-10}$   &  -9.301091200095 $\cdot 10^{-10}$   &  1.1 $\cdot 10^{-21}$ \\  \hline 
0.7   &   -6.989817680976 $\cdot 10^{-10}$   &  -6.987973346248 $\cdot 10^{-10}$   &  1.8 $\cdot 10^{-13}$ \\  \hline 
0.8   &   -4.652100850091 $\cdot 10^{-10}$   &  -4.652100850080 $\cdot 10^{-10}$   &  1.1 $\cdot 10^{-21}$ \\  \hline 
0.9   &   -2.319273360260 $\cdot 10^{-10}$   &  -2.321569863571 $\cdot 10^{-10}$   &  2.2 $\cdot 10^{-13}$ \\  \hline 
1   &     0   &   0                                 &  0 \\  \hline
    \end{tabular}}
\end{center}
\end{table}
\clearpage
\begin{table}[h!]\label{table7}
\textbf{Table 7.} Comparison between the HPM results
\cite{Moghimi2011}, OHPM results (\ref{JefferyHamel56}) and
numerical results for the velocity $F(\eta)$ for $\alpha =
\frac{\pi}{24}$ and $H=1000$ \\
\begin{center}
\scriptsize{
\begin{tabular}{|c|c|c|c|c|}
  \hline
  $\eta$ & $F_{\textrm{HPM}}$ \cite{Moghimi2011} &   $F_{\textrm{numeric}}$ & $ \bar{F}_{\textrm{OHPM}}$, Eq. (\ref{JefferyHamel56})& $\begin{array}{c}
                                                         \textrm{relative error}= \\
                                                          |F_{\textrm{numeric}}-\bar{F}_{\textrm{OHPM}}|
                                                          \end{array}$\\
                                                            \hline
  0   & 1  &   1   &   1   &   0 \\  \hline 
0.1   & 0.9937607  &   0.9937611413   &   0.9937578349   &   3.3 $\cdot 10^{-6}$ \\  \hline 
0.2   & 0.9747886  &   0.9747905377   &   0.9747871858   &   3.3 $\cdot 10^{-6}$ \\  \hline 
0.3   & 0.9422794  &   0.9422838643   &   0.9422837661   &   9.8 $\cdot 10^{-8}$ \\  \hline 
0.4   & 0.8947431  &   0.8947513541   &   0.8947499642   &   1.3 $\cdot 10^{-6}$ \\  \hline 
0.5   & 0.8297471  &   0.8297605618   &   0.8297564399   &   4.1 $\cdot 10^{-6}$ \\  \hline 
0.6   & 0.7434756  &   0.7434959928   &   0.7434941281   &   1.8 $\cdot 10^{-6}$ \\  \hline 
0.7   & 0.6299866  &   0.6300164190   &   0.6300167936   &   3.7 $\cdot 10^{-7}$ \\  \hline 
0.8   & 0.4799338  &   0.4799757286   &   0.4799724842   &   3.2 $\cdot 10^{-6}$ \\  \hline 
0.9   & 0.2782789  &   0.2783297889   &   0.2783280591   &   1.7 $\cdot 10^{-6}$ \\  \hline 
1    &  0  &   0  &    0   &    0 \\  \hline
    \end{tabular}}
\end{center}
\end{table}
\begin{table}[h!]\label{table8}
\textbf{Table 8.} Comparison between OHPM results
(\ref{JefferyHamel57}) and numerical results for the
temperature $\theta(\eta)$ for $\alpha = \frac{\pi}{24}$ and $H=1000$ \\
\begin{center}
\scriptsize{
\begin{tabular}{|c|c|c|c|}
  \hline
  $\eta$ & ${\theta}_{\textrm{numeric}}$ & $ {\bar{{\theta}}}_{\textrm{OHPM}}$ from Eq. (\ref{JefferyHamel57})& $\begin{array}{c}
                                                         \textrm{relative error}= \\
                                                          |{{\theta}}_{\textrm{numeric}}-{\bar{{\theta}}}_{\textrm{OHPM}}|
                                                          \end{array}$\\
                                                            \hline
  0   &   -3.885903481818 $\cdot 10^{-9}$   &   -3.885903481801 $\cdot 10^{-9}$   &   1.7 $\cdot 10^{-20}$ \\  \hline 
0.1   &   -3.804390715058 $\cdot 10^{-9}$   &   -3.805365047292 $\cdot 10^{-9}$   &   9.7 $\cdot 10^{-13}$ \\  \hline 
0.2   &   -3.575104095175 $\cdot 10^{-9}$   &   -3.575104095175 $\cdot 10^{-9}$   &   5.7 $\cdot 10^{-23}$ \\  \hline 
0.3   &   -3.222962704714 $\cdot 10^{-9}$   &   -3.222630117770 $\cdot 10^{-9}$   &   3.3 $\cdot 10^{-13}$ \\  \hline 
0.4   &   -2.782964569511 $\cdot 10^{-9}$   &   -2.782964569511 $\cdot 10^{-9}$   &   4.4 $\cdot 10^{-23}$ \\  \hline 
0.5   &   -2.293213461698 $\cdot 10^{-9}$   &   -2.293504856351 $\cdot 10^{-9}$   &   2.9 $\cdot 10^{-13}$ \\  \hline 
0.6   &   -1.790055410191 $\cdot 10^{-9}$   &   -1.790055410191 $\cdot 10^{-9}$   &   3.5 $\cdot 10^{-23}$ \\  \hline 
0.7   &   -1.302936492610 $\cdot 10^{-9}$   &   -1.301687731203 $\cdot 10^{-9}$   &   1.2 $\cdot 10^{-12}$ \\  \hline 
0.8   &   -8.442618285524 $\cdot 10^{-10}$   &  -8.442618285525 $\cdot 10^{-10}$   &  1.4 $\cdot 10^{-23}$ \\  \hline 
0.9   &   -4.107432220830 $\cdot 10^{-10}$   &  -4.158108458702 $\cdot 10^{-10}$   &  5.06 $\cdot 10^{-12}$ \\  \hline 
1    &     0   &   0                                 &  0 \\  \hline
    \end{tabular}}
\end{center}
\end{table}
\clearpage
\begin{table}[h!]\label{table9}
\textbf{Table 9.} Comparison between the OHPM results
(\ref{JefferyHamel58}) and numerical results for the velocity
$F(\eta)$ for $\alpha =
\frac{\pi}{36}$ and $H=0$ \\
\begin{center}
\scriptsize{
\begin{tabular}{|c|c|c|c|}
  \hline
  $\eta$ & $F_{\textrm{numeric}}$ & $ \bar{F}_{\textrm{OHPM}}$, Eq. (\ref{JefferyHamel58})& $\begin{array}{c}
                                                         \textrm{relative error}= \\
                                                          |F_{\textrm{numeric}}-\bar{F}_{\textrm{OHPM}}|
                                                          \end{array}$\\
                                                            \hline
  0   &   1   &   1   &   0 \\  \hline  
0.1   &   0.9824312364   &   0.9824320701   &   8.3 $\cdot 10^{-7}$ \\  \hline  
0.2   &   0.9312259577   &   0.9312265466   &   5.8 $\cdot 10^{-7}$ \\  \hline  
0.3   &   0.8506106161   &   0.8506107339   &   1.1 $\cdot 10^{-7}$ \\  \hline  
0.4   &   0.7467908018   &   0.7467915008   &   6.9 $\cdot 10^{-7}$ \\  \hline  
0.5   &   0.6269481682   &   0.6269490072   &   8.3 $\cdot 10^{-7}$ \\  \hline  
0.6   &   0.4982344464   &   0.4982347463   &   2.9 $\cdot 10^{-7}$ \\  \hline  
0.7   &   0.3669663386   &   0.3669668037   &   4.6 $\cdot 10^{-7}$ \\  \hline  
0.8   &   0.2381237463   &   0.2381245686   &   8.2 $\cdot 10^{-7}$ \\  \hline  
0.9   &   0.1151519312   &   0.1151523953   &   4.6 $\cdot 10^{-7}$ \\  \hline  
  1   &   0   &   0   &     0  \\  \hline  
    \end{tabular}}
\end{center}
\end{table}
\begin{table}[h!]\label{table10}
\textbf{Table 10.} Comparison between OHPM results
(\ref{JefferyHamel59}) and numerical results for the
temperature $\theta(\eta)$ for $\alpha = \frac{\pi}{36}$ and $H=0$ \\
\begin{center}
\scriptsize{
\begin{tabular}{|c|c|c|c|}
  \hline
  $\eta$ & ${\theta}_{\textrm{numeric}}$ & $ {\bar{{\theta}}}_{\textrm{OHPM}}$ from Eq. (\ref{JefferyHamel59})& $\begin{array}{c}
                                                         \textrm{relative error}= \\
                                                          |{{\theta}}_{\textrm{numeric}}-{\bar{{\theta}}}_{\textrm{OHPM}}|
                                                          \end{array}$\\
                                                            \hline
  0   &   -2.552530215568 $\cdot 10^{-11}$   &   -2.552530214568 $\cdot 10^{-11}$   &   9.9 $\cdot 10^{-21}$ \\  \hline  
0.1   &   -2.442980325655 $\cdot 10^{-11}$   &   -2.444079060328 $\cdot 10^{-11}$   &   1.09 $\cdot 10^{-14}$ \\  \hline  
0.2   &   -2.143998685224 $\cdot 10^{-11}$   &   -2.143998685224 $\cdot 10^{-11}$   &   1.9 $\cdot 10^{-25}$ \\  \hline  
0.3   &   -1.714518106214 $\cdot 10^{-11}$   &   -1.713385692946 $\cdot 10^{-11}$   &   1.1 $\cdot 10^{-14}$ \\  \hline  
0.4   &   -1.227729149275 $\cdot 10^{-11}$   &   -1.227729149275 $\cdot 10^{-11}$   &   1.3 $\cdot 10^{-25}$ \\  \hline  
0.5   &   -7.579058130819 $\cdot 10^{-12}$   &   -7.588196380463 $\cdot 10^{-12}$   &   9.1 $\cdot 10^{-15}$ \\  \hline  
0.6   &   -3.636709774083 $\cdot 10^{-12}$   &   -3.636709774083 $\cdot 10^{-12}$   &   1.6 $\cdot 10^{-25}$ \\  \hline  
0.7   &   -8.316369531201 $\cdot 10^{-13}$   &   -8.099883012752 $\cdot 10^{-13}$   &   2.1 $\cdot 10^{-14}$ \\  \hline  
0.8   &    6.930669117891 $\cdot 10^{-13}$   &    6.930669117890 $\cdot 10^{-13}$   &   7.8 $\cdot 10^{-26}$ \\  \hline  
0.9   &    9.721365514984 $\cdot 10^{-13}$   &    8.898023093848 $\cdot 10^{-13}$   &   8.2 $\cdot 10^{-14}$ \\  \hline  
  1   &    0  &  0    &   0  \\  \hline
    \end{tabular}}
\end{center}
\end{table}
\clearpage
\begin{table}[h!]\label{table11}
\textbf{Table 11.} Comparison between the OHPM results
(\ref{JefferyHamel60}) and numerical results for the velocity
$F(\eta)$ for $\alpha =
\frac{\pi}{36}$ and $H=250$ \\
\begin{center}
\scriptsize{
\begin{tabular}{|c|c|c|c|}
  \hline
  $\eta$ &  $F_{\textrm{numeric}}$ & $ \bar{F}_{\textrm{OHPM}}$, Eq. (\ref{JefferyHamel60})& $\begin{array}{c}
                                                         \textrm{relative error}= \\
                                                          |F_{\textrm{numeric}}-\bar{F}_{\textrm{OHPM}}|
                                                          \end{array}$\\
                                                            \hline
  0   &   1   &   1   &   0 \\  \hline  
0.1   &   0.9849746347   &   0.9849746827   &   4.8 $\cdot 10^{-8}$ \\  \hline  
0.2   &   0.9409030596   &   0.9409030371   &   2.2 $\cdot 10^{-8}$ \\  \hline  
0.3   &   0.8706210985   &   0.8706210491   &   4.9 $\cdot 10^{-8}$ \\  \hline  
0.4   &   0.7783094304   &   0.7783094038   &   2.6 $\cdot 10^{-8}$ \\  \hline  
0.5   &   0.6688194854   &   0.6688193877   &   9.7 $\cdot 10^{-8}$ \\  \hline  
0.6   &   0.5469607278   &   0.5469605971   &   1.3 $\cdot 10^{-7}$ \\  \hline  
0.7   &   0.4168757453   &   0.4168756336   &   1.1 $\cdot 10^{-7}$ \\  \hline  
0.8   &   0.2815737222   &   0.2815735380   &   1.8 $\cdot 10^{-7}$ \\  \hline  
0.9   &   0.1426291045   &   0.1426289244   &   1.8 $\cdot 10^{-7}$ \\  \hline  
  1   &   0    &     0     &   0   \\  \hline 
    \end{tabular}}
\end{center}
\end{table}
\begin{table}[h!]\label{table12}
\textbf{Table 12.} Comparison between OHPM results
(\ref{JefferyHamel61}) and numerical results for the
temperature $\theta(\eta)$ for $\alpha = \frac{\pi}{36}$ and $H=250$ \\
\begin{center}
\scriptsize{
\begin{tabular}{|c|c|c|c|}
  \hline
  $\eta$ & ${\theta}_{\textrm{numeric}}$ & $ {\bar{{\theta}}}_{\textrm{OHPM}}$ from Eq. (\ref{JefferyHamel61})& $\begin{array}{c}
                                                         \textrm{relative error}= \\
                                                          |{{\theta}}_{\textrm{numeric}}-{\bar{{\theta}}}_{\textrm{OHPM}}|
                                                          \end{array}$\\
                                                            \hline
  0   &   -3.397742006028 $\cdot 10^{-9}$   &   -3.397742006018 $\cdot 10^{-9}$   &   9.9 $\cdot 10^{-21}$ \\  \hline  
0.1   &   -3.346638363267 $\cdot 10^{-9}$   &   -3.347192325339 $\cdot 10^{-9}$   &   5.5 $\cdot 10^{-13}$ \\  \hline  
0.2   &   -3.199501303768 $\cdot 10^{-9}$   &   -3.199501303768 $\cdot 10^{-9}$   &   1.9 $\cdot 10^{-22}$ \\  \hline  
0.3   &   -2.964278501171 $\cdot 10^{-9}$   &   -2.964072110215 $\cdot 10^{-9}$   &   2.06 $\cdot 10^{-13}$ \\  \hline  
0.4   &   -2.653181572584 $\cdot 10^{-9}$   &   -2.653181572584 $\cdot 10^{-9}$   &   6.7 $\cdot 10^{-23}$ \\  \hline  
0.5   &   -2.281170448213 $\cdot 10^{-9}$   &   -2.281224162207 $\cdot 10^{-9}$   &   5.3 $\cdot 10^{-14}$ \\  \hline  
0.6   &   -1.864037686164 $\cdot 10^{-9}$   &   -1.864037686164 $\cdot 10^{-9}$   &   1.4 $\cdot 10^{-23}$ \\  \hline  
0.7   &   -1.416766717489 $\cdot 10^{-9}$   &   -1.416988166987 $\cdot 10^{-9}$   &   2.2 $\cdot 10^{-13}$ \\  \hline  
0.8   &   -9.521674889800 $\cdot 10^{-10}$   &  -9.521674889800 $\cdot 10^{-10}$   &  2.8 $\cdot 10^{-24}$ \\  \hline  
0.9   &   -4.798620033623 $\cdot 10^{-10}$   &  -4.772452449206 $\cdot 10^{-10}$   &  2.6 $\cdot 10^{-12}$ \\  \hline  
  1   &   0   &   0   &   0  \\  \hline
    \end{tabular}}
\end{center}
\end{table}
\clearpage
\begin{table}[h!]\label{table13}
\textbf{Table 13.} Comparison between the OHPM results
(\ref{JefferyHamel62}) and numerical results for the velocity
$F(\eta)$ for $\alpha =
\frac{\pi}{36}$ and $H=500$ \\
\begin{center}
\scriptsize{
\begin{tabular}{|c|c|c|c|}
  \hline
  $\eta$ &  $F_{\textrm{numeric}}$ & $ \bar{F}_{\textrm{OHPM}}$, Eq. (\ref{JefferyHamel62})& $\begin{array}{c}
                                                         \textrm{relative error}= \\
                                                          |F_{\textrm{numeric}}-\bar{F}_{\textrm{OHPM}}|
                                                          \end{array}$\\
                                                            \hline
  0   &   1   &   1   &   0 \\  \hline  
0.1   &   0.9871108049   &   0.9871108182   &   1.3 $\cdot 10^{-8}$ \\  \hline  
0.2   &   0.9490642316   &   0.9490642266   &   5.0 $\cdot 10^{-9}$ \\  \hline  
0.3   &   0.8876104540   &   0.8876104486   &   5.3 $\cdot 10^{-9}$ \\  \hline  
0.4   &   0.8053144512   &   0.8053144616   &   1.03 $\cdot 10^{-8}$ \\  \hline  
0.5   &   0.7051032297   &   0.7051032241   &   5.5 $\cdot 10^{-9}$ \\  \hline  
0.6   &   0.5897534597   &   0.5897534552   &   4.4 $\cdot 10^{-9}$ \\  \hline  
0.7   &   0.4613867287   &   0.4613867398   &   1.1 $\cdot 10^{-8}$ \\  \hline  
0.8   &   0.3210064088   &   0.3210063961   &   1.2 $\cdot 10^{-8}$ \\  \hline  
0.9   &   0.1680649652   &   0.1680649814   &   1.6 $\cdot 10^{-8}$ \\  \hline  
 1   &   0   &  0   &   0   \\  \hline 
    \end{tabular}}
\end{center}
\end{table}
\begin{table}[h!]\label{table14}
\textbf{Table 14.} Comparison between OHPM results
(\ref{JefferyHamel63}) and numerical results for the
temperature $\theta(\eta)$ for $\alpha = \frac{\pi}{36}$ and $H=500$ \\
\begin{center}
\scriptsize{
\begin{tabular}{|c|c|c|c|}
  \hline
  $\eta$ & ${\theta}_{\textrm{numeric}}$ & $ {\bar{{\theta}}}_{\textrm{OHPM}}$ from Eq. (\ref{JefferyHamel63})& $\begin{array}{c}
                                                         \textrm{relative error}= \\
                                                          |{{\theta}}_{\textrm{numeric}}-{\bar{{\theta}}}_{\textrm{OHPM}}|
                                                          \end{array}$\\
                                                            \hline
  0   &   -6.861779213872 $\cdot 10^{-9}$   &   -6.861779213862 $\cdot 10^{-9}$   &   9.9 $\cdot 10^{-21}$ \\  \hline  
0.1   &   -6.756689795935 $\cdot 10^{-9}$   &   -6.757837516749 $\cdot 10^{-9}$   &   1.1 $\cdot 10^{-12}$ \\  \hline  
0.2   &   -6.454596023395 $\cdot 10^{-9}$   &   -6.454596023395 $\cdot 10^{-9}$   &   1.4 $\cdot 10^{-23}$ \\  \hline  
0.3   &   -5.973076729675 $\cdot 10^{-9}$   &   -5.972618579646 $\cdot 10^{-9}$   &   4.5 $\cdot 10^{-13}$ \\  \hline  
0.4   &   -5.338639377096 $\cdot 10^{-9}$   &   -5.338639377096 $\cdot 10^{-9}$   &   5.7 $\cdot 10^{-24}$ \\  \hline  
0.5   &   -4.583265280198 $\cdot 10^{-9}$   &   -4.583356882876 $\cdot 10^{-9}$   &   9.1 $\cdot 10^{-14}$ \\  \hline  
0.6   &   -3.739726496008 $\cdot 10^{-9}$   &   -3.739726496008 $\cdot 10^{-9}$   &   6.6 $\cdot 10^{-24}$ \\  \hline  
0.7   &   -2.838920429095 $\cdot 10^{-9}$   &   -2.839180849233 $\cdot 10^{-9}$   &   2.6 $\cdot 10^{-13}$ \\  \hline  
0.8   &   -1.906149341683 $\cdot 10^{-9}$   &   -1.906149341683 $\cdot 10^{-9}$   &   5.7 $\cdot 10^{-24}$ \\  \hline  
0.9   &   -9.597066242185 $\cdot 10^{-10}$   &  -9.553857275377 $\cdot 10^{-10}$   &  4.3 $\cdot 10^{-12}$ \\  \hline  
  1   &   0   &   0     &    0       \\  \hline
    \end{tabular}}
\end{center}
\end{table}
\clearpage
\begin{table}[h!]\label{table15}
\textbf{Table 15.} Comparison between the OHPM results
(\ref{JefferyHamel64}) and numerical results for the velocity
$F(\eta)$ for $\alpha =
\frac{\pi}{36}$ and $H=1000$ \\
\begin{center}
\scriptsize{
\begin{tabular}{|c|c|c|c|}
  \hline
  $\eta$ &  $F_{\textrm{numeric}}$ & $ \bar{F}_{\textrm{OHPM}}$, Eq. (\ref{JefferyHamel64})& $\begin{array}{c}
                                                         \textrm{relative error}= \\
                                                          |F_{\textrm{numeric}}-\bar{F}_{\textrm{OHPM}}|
                                                          \end{array}$\\
                                                            \hline
  0   &   1   &   1   &   0 \\  \hline  
0.1   &   0.9904272110   &   0.9904271107   &   1.003 $\cdot 10^{-7}$ \\  \hline  
0.2   &   0.9618100677   &   0.9618099299   &   1.3 $\cdot 10^{-7}$ \\  \hline  
0.3   &   0.9144036356   &   0.9144036639   &   2.8 $\cdot 10^{-8}$ \\  \hline  
0.4   &   0.8484736891   &   0.8484737167   &   2.7 $\cdot 10^{-8}$ \\  \hline  
0.5   &   0.7640642211   &   0.7640640849   &   1.3 $\cdot 10^{-7}$ \\  \hline  
0.6   &   0.6606772356   &   0.6606771839   &   5.1 $\cdot 10^{-8}$ \\  \hline  
0.7   &   0.5368520578   &   0.5368521821   &   1.2 $\cdot 10^{-7}$ \\  \hline  
0.8   &   0.3896018441   &   0.3896017829   &   6.1 $\cdot 10^{-8}$ \\  \hline  
0.9   &   0.2136111325   &   0.2136111294   &   3.1 $\cdot 10^{-9}$ \\  \hline  
 1   &   0   &   0    &   0   \\  \hline 
    \end{tabular}}
\end{center}
\end{table}
\begin{table}[h!]\label{table16}
\textbf{Table 16.} Comparison between OHPM results
(\ref{JefferyHamel65}) and numerical results for the
temperature $\theta(\eta)$ for $\alpha = \frac{\pi}{36}$ and $H=1000$ \\
\begin{center}
\scriptsize{
\begin{tabular}{|c|c|c|c|}
  \hline
  $\eta$ & ${\theta}_{\textrm{numeric}}$ & $ {\bar{{\theta}}}_{\textrm{OHPM}}$ from Eq. (\ref{JefferyHamel65})& $\begin{array}{c}
                                                         \textrm{relative error}= \\
                                                          |{{\theta}}_{\textrm{numeric}}-{\bar{{\theta}}}_{\textrm{OHPM}}|
                                                          \end{array}$\\
                                                            \hline
  0   &   -1.406496797937 $\cdot 10^{-8}$   &   -1.406496797936 $\cdot 10^{-8}$   &   9.9 $\cdot 10^{-21}$ \\  \hline  
0.1   &   -1.383991714524 $\cdot 10^{-8}$   &   -1.384237616580 $\cdot 10^{-8}$   &   2.4 $\cdot 10^{-12}$ \\  \hline  
0.2   &   -1.319483359510 $\cdot 10^{-8}$   &   -1.319483359510 $\cdot 10^{-8}$   &   8.2 $\cdot 10^{-24}$ \\  \hline  
0.3   &   -1.217244871718 $\cdot 10^{-8}$   &   -1.217139646331 $\cdot 10^{-8}$   &   1.05 $\cdot 10^{-12}$ \\  \hline  
0.4   &   -1.083591345359 $\cdot 10^{-8}$   &   -1.083591345359 $\cdot 10^{-8}$   &   1.1 $\cdot 10^{-23}$ \\  \hline  
0.5   &   -9.260353523538 $\cdot 10^{-9}$   &   -9.260364501045 $\cdot 10^{-9}$   &   1.09 $\cdot 10^{-14}$ \\  \hline  
0.6   &   -7.519751014022 $\cdot 10^{-9}$   &   -7.519751014022 $\cdot 10^{-9}$   &   4.9 $\cdot 10^{-24}$ \\  \hline  
0.7   &   -5.683228753360 $\cdot 10^{-9}$   &   -5.683311365697 $\cdot 10^{-9}$   &   8.2 $\cdot 10^{-14}$ \\  \hline  
0.8   &   -3.802282135607 $\cdot 10^{-9}$   &   -3.802282135607 $\cdot 10^{-9}$   &   3.3 $\cdot 10^{-24}$ \\  \hline  
0.9   &   -1.908111552553 $\cdot 10^{-9}$   &   -1.902733734257 $\cdot 10^{-9}$   &   5.3 $\cdot 10^{-12}$ \\  \hline  
  1   &   0   &    0    &    0     \\  \hline
    \end{tabular}}
\end{center}
\end{table}
\clearpage
\noindent
\begin{tabular}[!t]{c}
\begin{tabular}[!t]{c c}
  \includegraphics[width=2.55in]{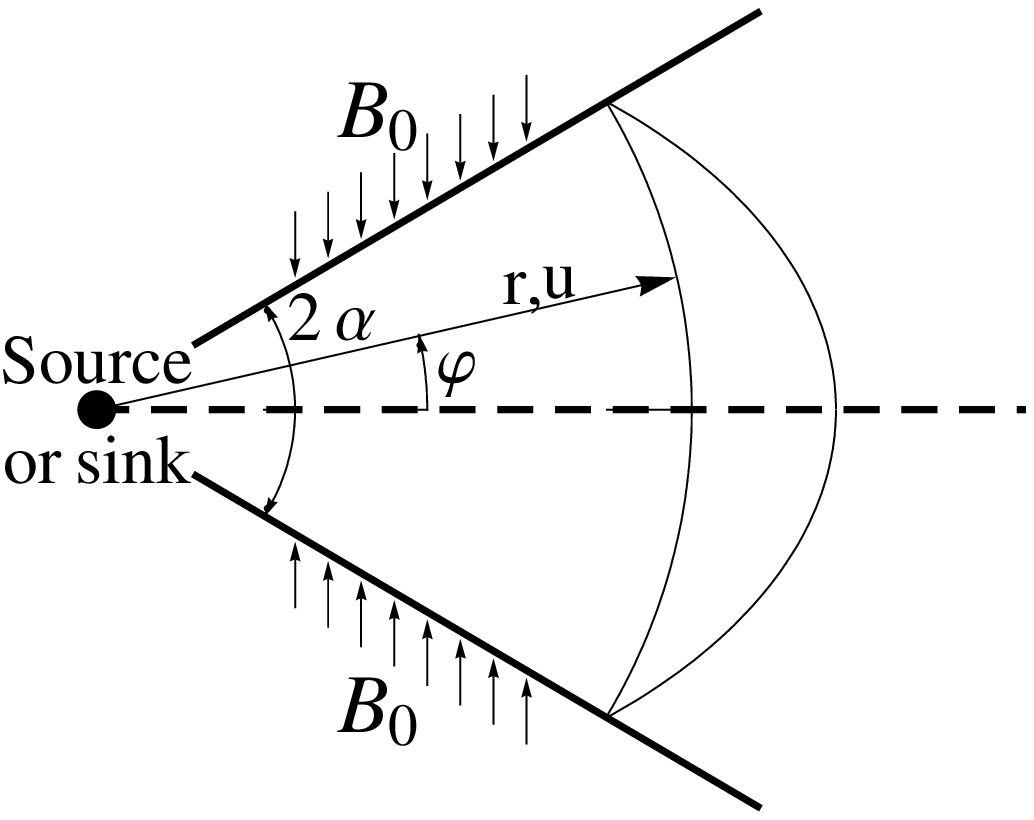} \\
\end{tabular}\\
\centerline {Fig 1.  Geometry of the MHD Jeffery-Hamel problem.}
\end{tabular}
\vspace{0.5 cm} \noindent
\begin{tabular}[!t]{c}
\begin{tabular}[!t]{c c}
  \includegraphics[width=2.15in]{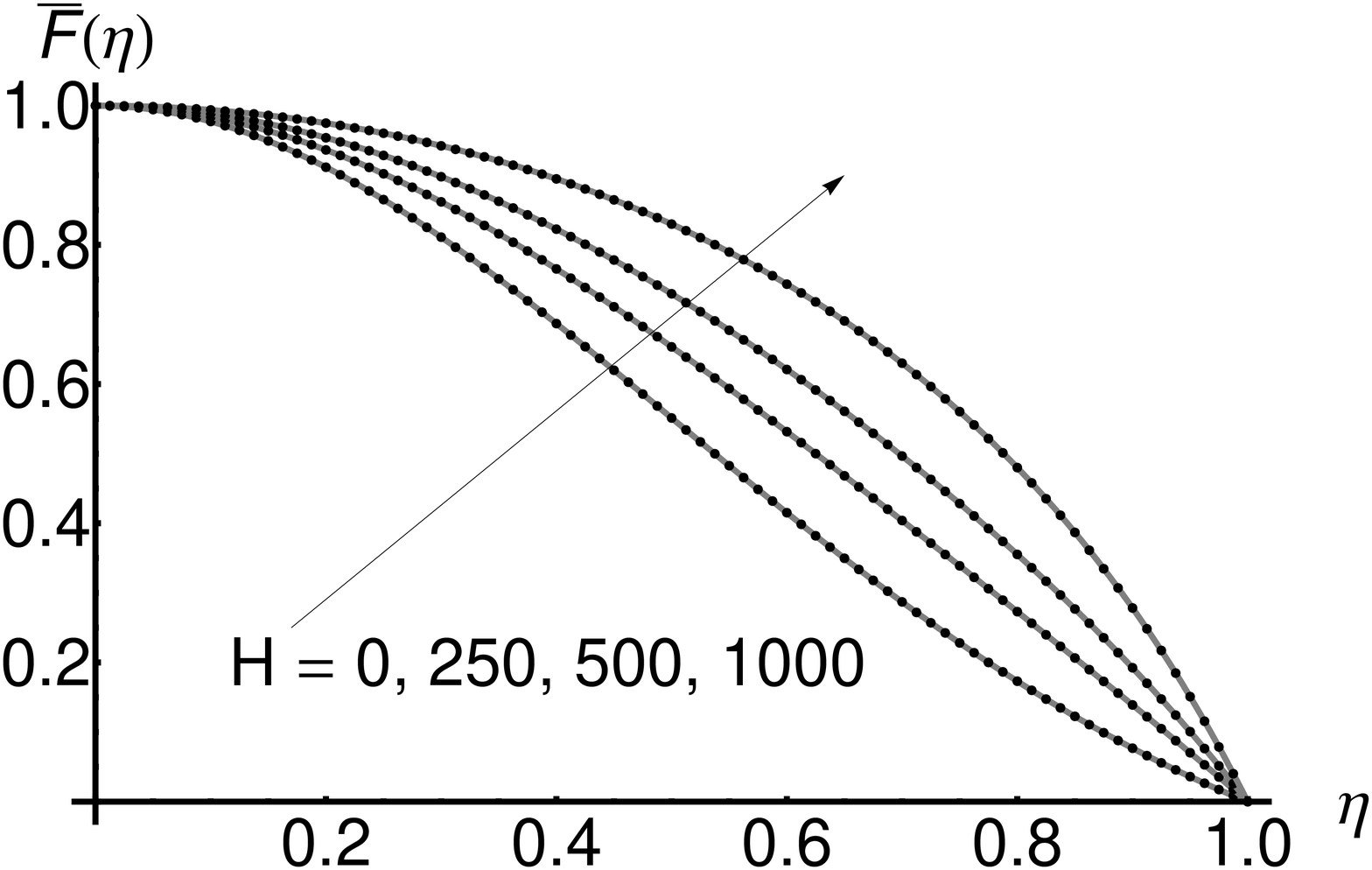}& \includegraphics[width=2.15in]{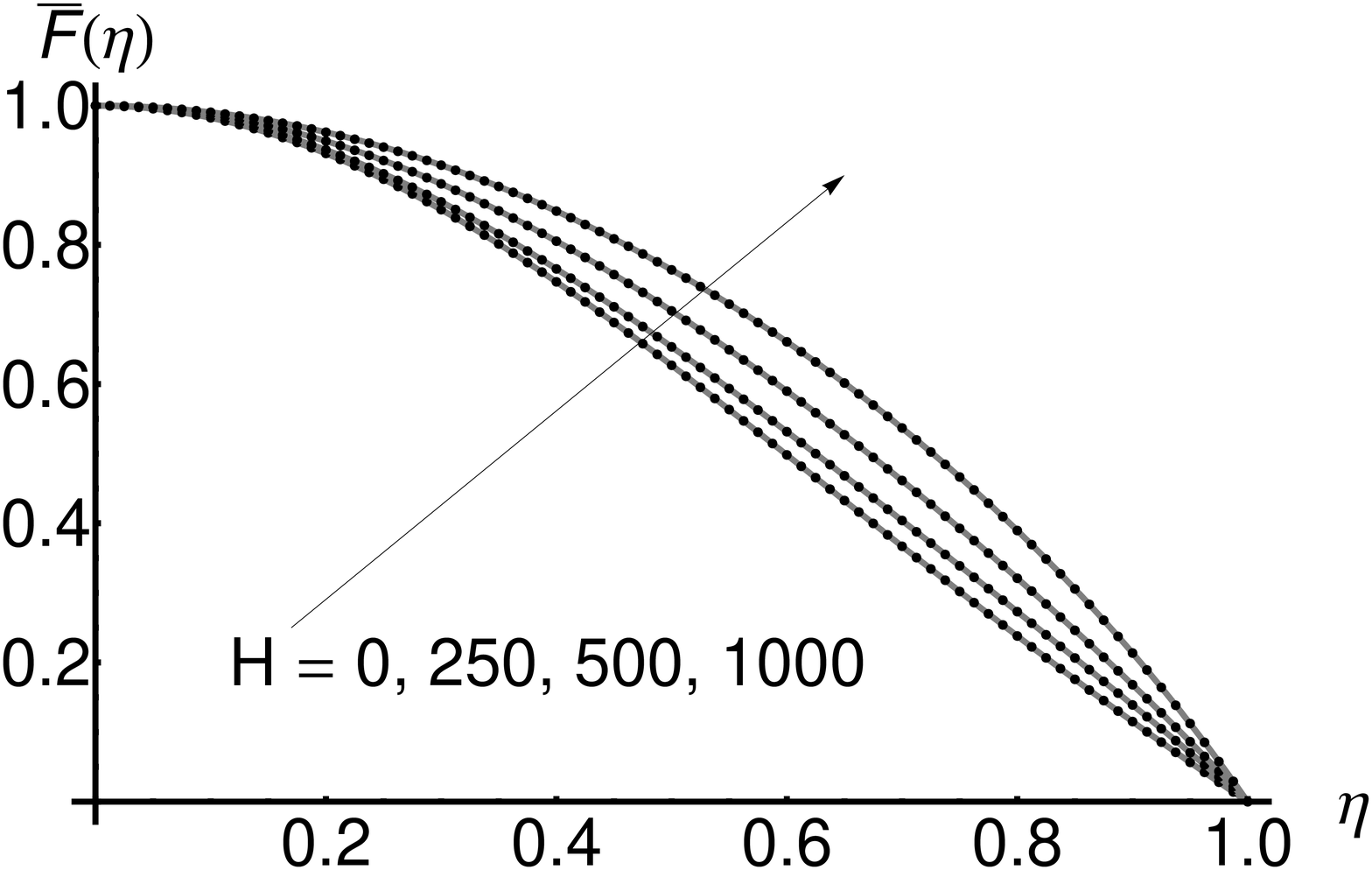} \\
{Fig. 2 Effect of the Hartmann }     &    Fig. 3 Effect of the Hartmann  \\
number on the velocity profile for      &      number  on the velocity profile for \\
$\alpha = \pi/24$, $Re=50$:   &   $\alpha = \pi/36$, $Re=50$: \\
------  numerical solution,      &    ------  numerical solution,       \\
\end{tabular}\\
\end{tabular}\\
\vspace{0.5 cm} \noindent
\begin{tabular}[!t]{c}
\begin{tabular}[!t]{c c}
  \includegraphics[width=2.15in]{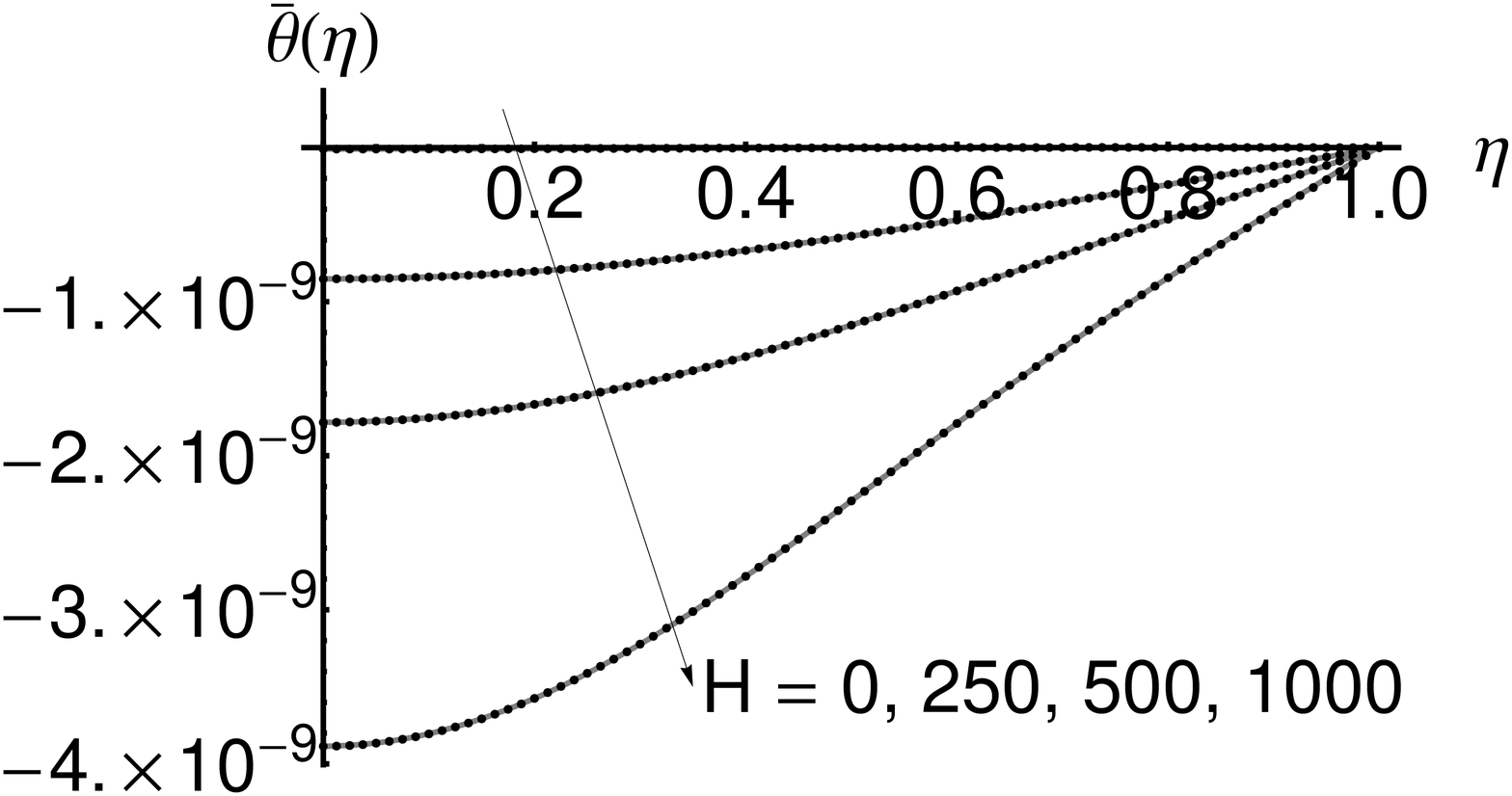}& \includegraphics[width=2.15in]{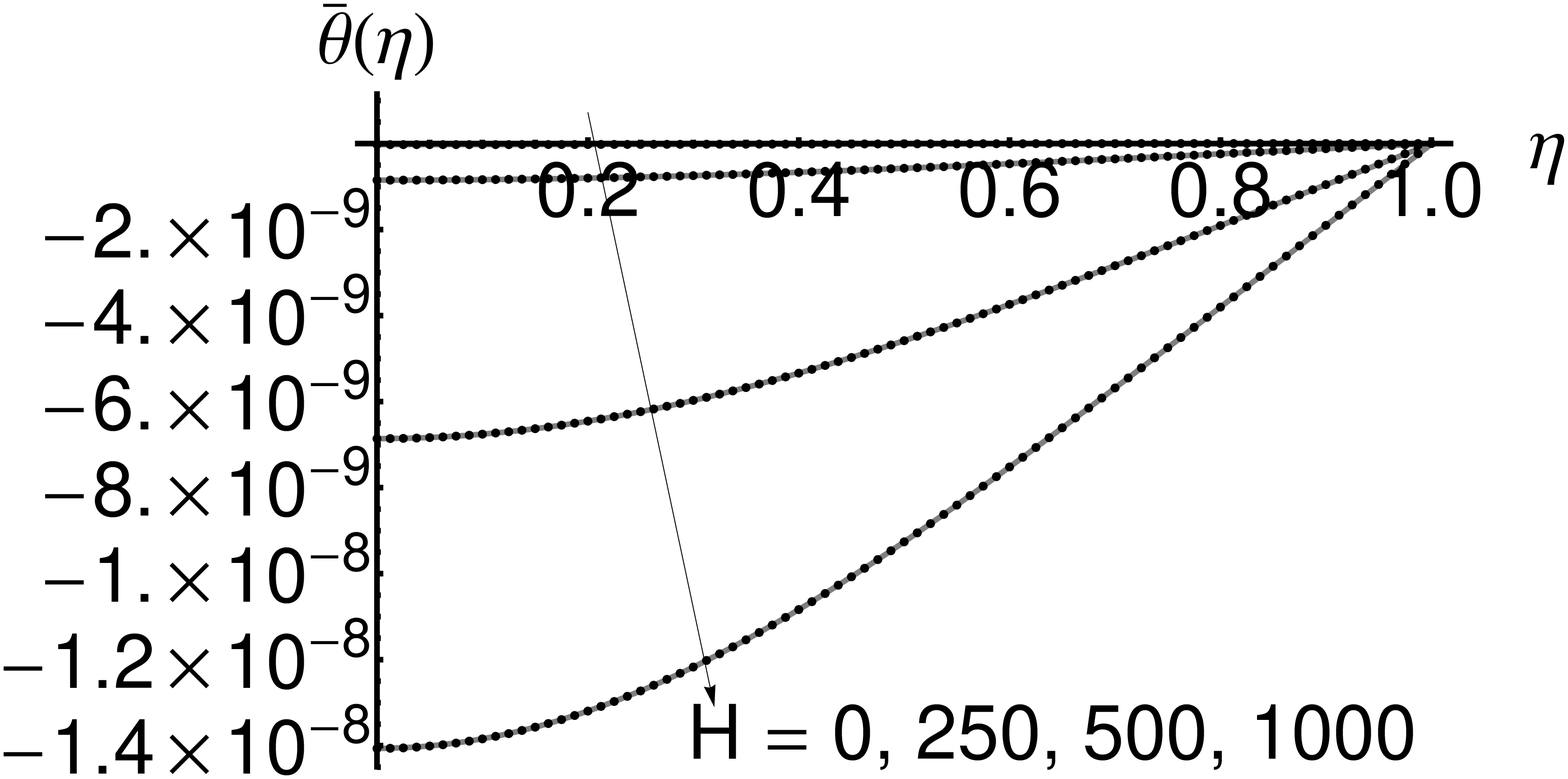} \\
Fig. 4 Effect of the Hartmann       &    Fig. 5 Effect of the Hartmann  \\
number on the thermal profile for      &      number  on the thermal profile for \\
$\alpha = \pi/24$, $Re=50$:   &   $\alpha = \pi/36$, $Re=50$: \\
------  numerical solution,      &    ------  numerical solution,       \\
........ OHPM solution &  
........ OHPM solution  \\
\end{tabular}\\
\end{tabular}\\
\begin{tabular}[!t]{c}
\begin{tabular}[!t]{c c}
  \includegraphics[width=2.15in]{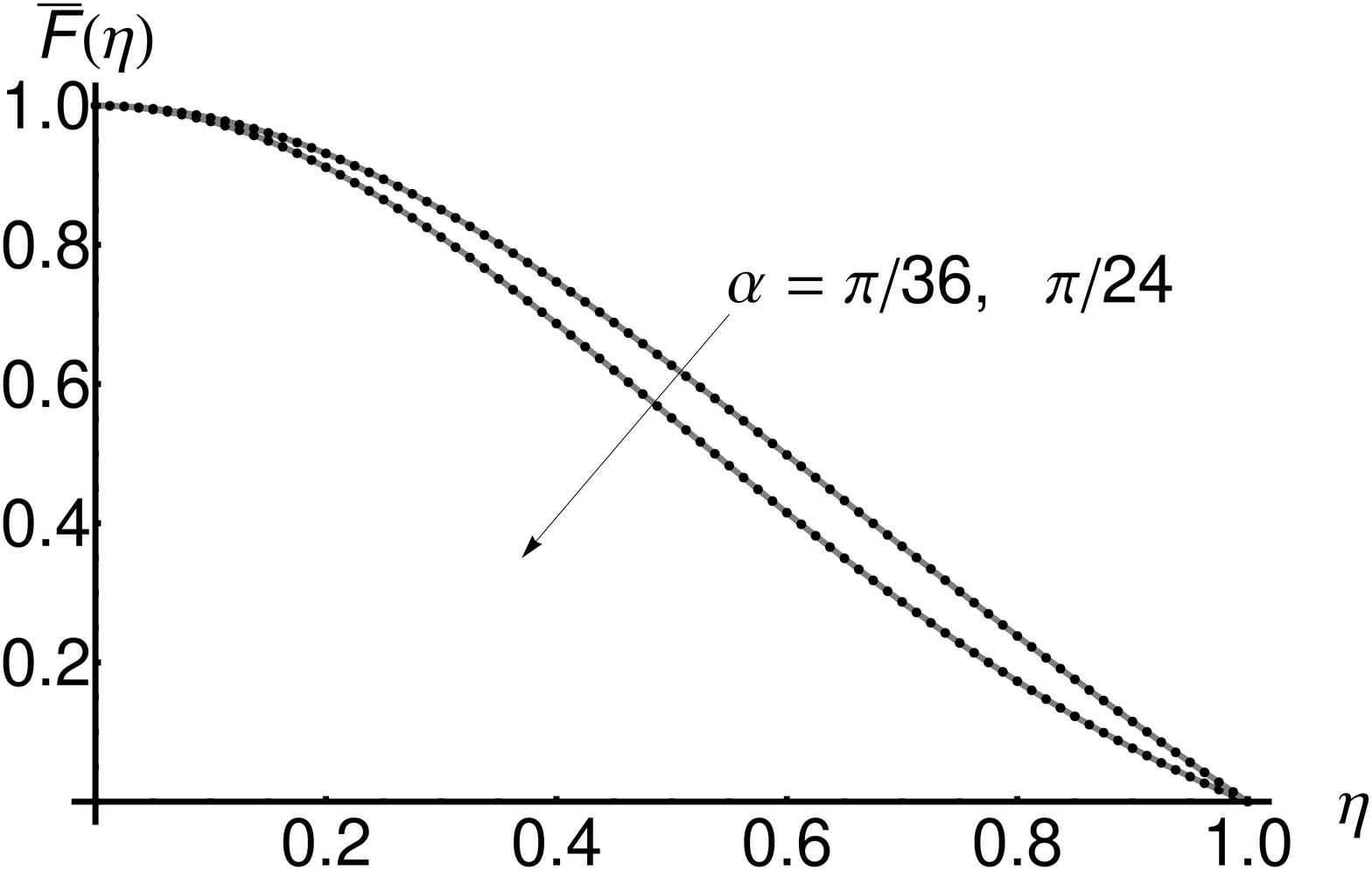}& \includegraphics[width=2.15in]{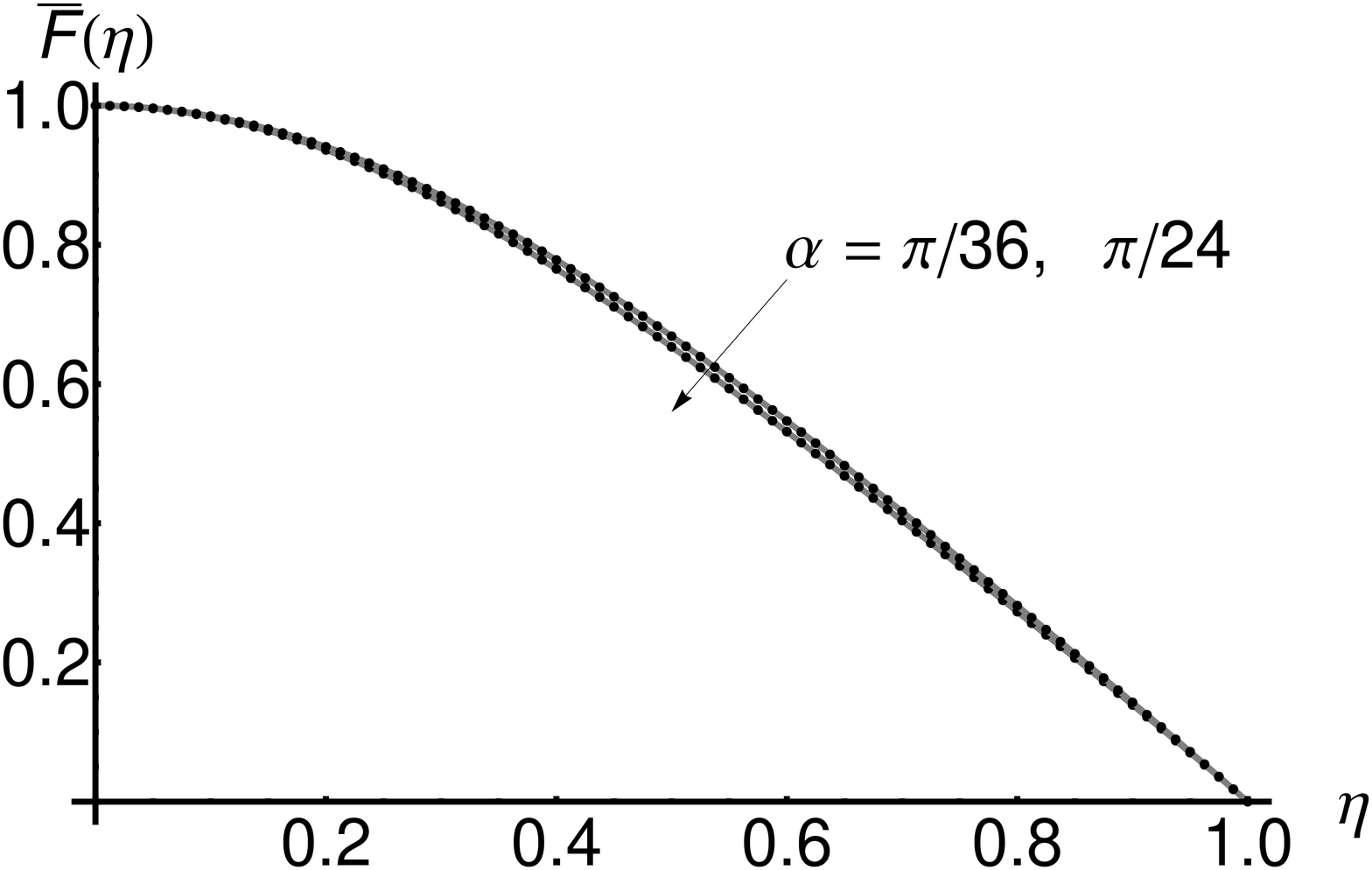}  \\
Fig. 6 Velocity profile  for $\alpha = \pi/36$,          &    Fig. 7 Velocity profile  for $\alpha = \pi/36$,     \\
and $\alpha = \pi/24$, $H=0$, $Re=50$,  &      and $\alpha = \pi/24$, $H=250$, $Re=50$, \\
Eqs. (\ref{JefferyHamel50}) and (\ref{JefferyHamel58}):   &   Eqs. (\ref{JefferyHamel52}) and (\ref{JefferyHamel60}): \\
------  numerical solution,      &    ------  numerical solution,       \\
........ OHPM solution &  
........ OHPM solution  \\
\end{tabular}\\
\end{tabular}
\begin{tabular}[!t]{c}
\begin{tabular}[!t]{c c}
  \includegraphics[width=2.15in]{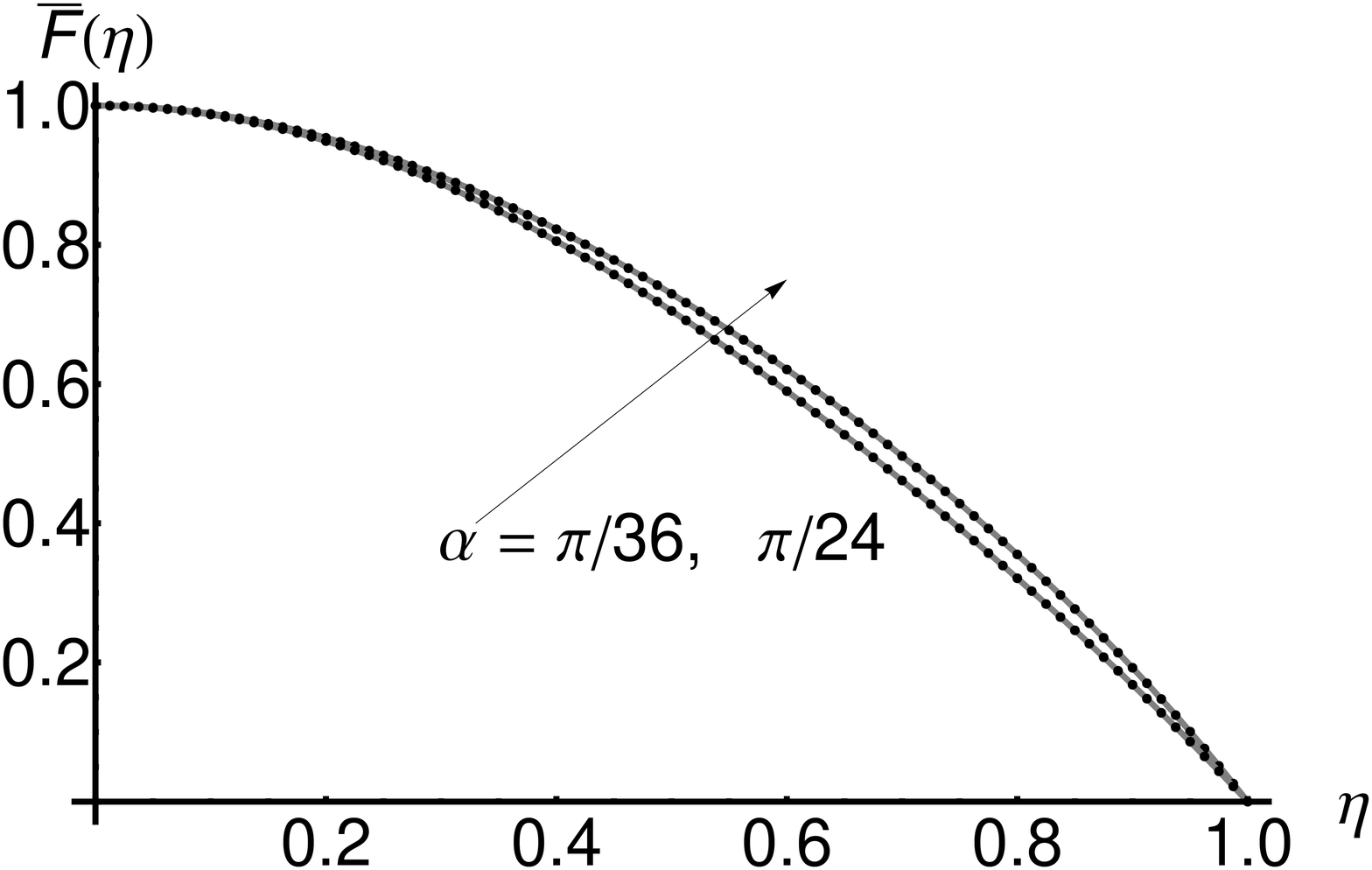}& \includegraphics[width=2.15in]{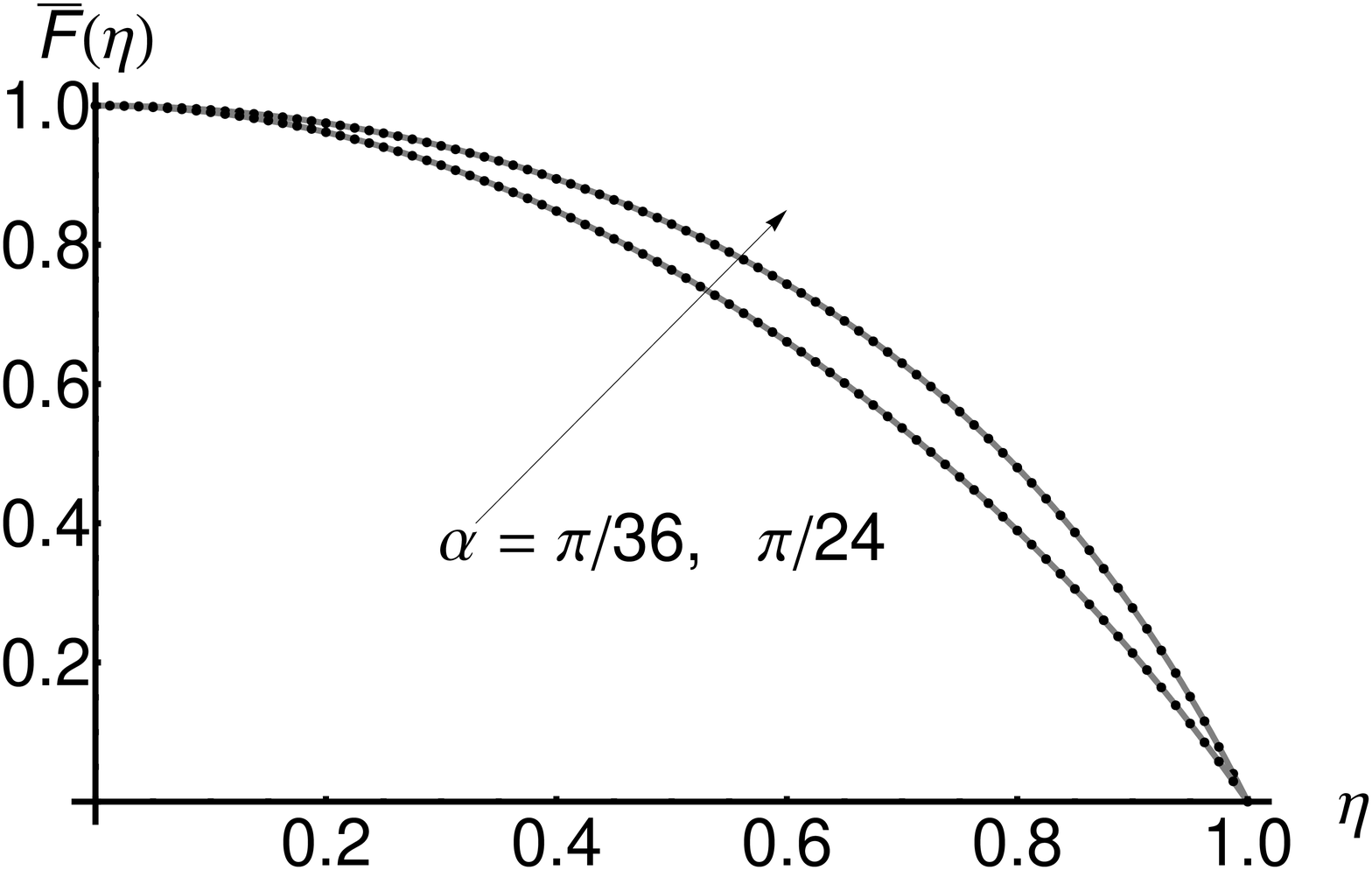}  \\
Fig. 8 Velocity profile  for $\alpha = \pi/36$,          &    Fig. 9 Velocity profile  for $\alpha = \pi/36$,     \\
and $\alpha = \pi/24$, $H=500$, $Re=50$,  &      and $\alpha = \pi/24$, $H=1000$, $Re=50$, \\
Eqs. (\ref{JefferyHamel54}) and (\ref{JefferyHamel62}):   &   Eqs. (\ref{JefferyHamel56}) and (\ref{JefferyHamel64}): \\
------  numerical solution,      &    ------  numerical solution,       \\
........ OHPM solution &  
........ OHPM solution  \\
\end{tabular}\\
\end{tabular}\\
\begin{tabular}[!t]{c}
\begin{tabular}[!t]{c c}
  \includegraphics[width=2.15in]{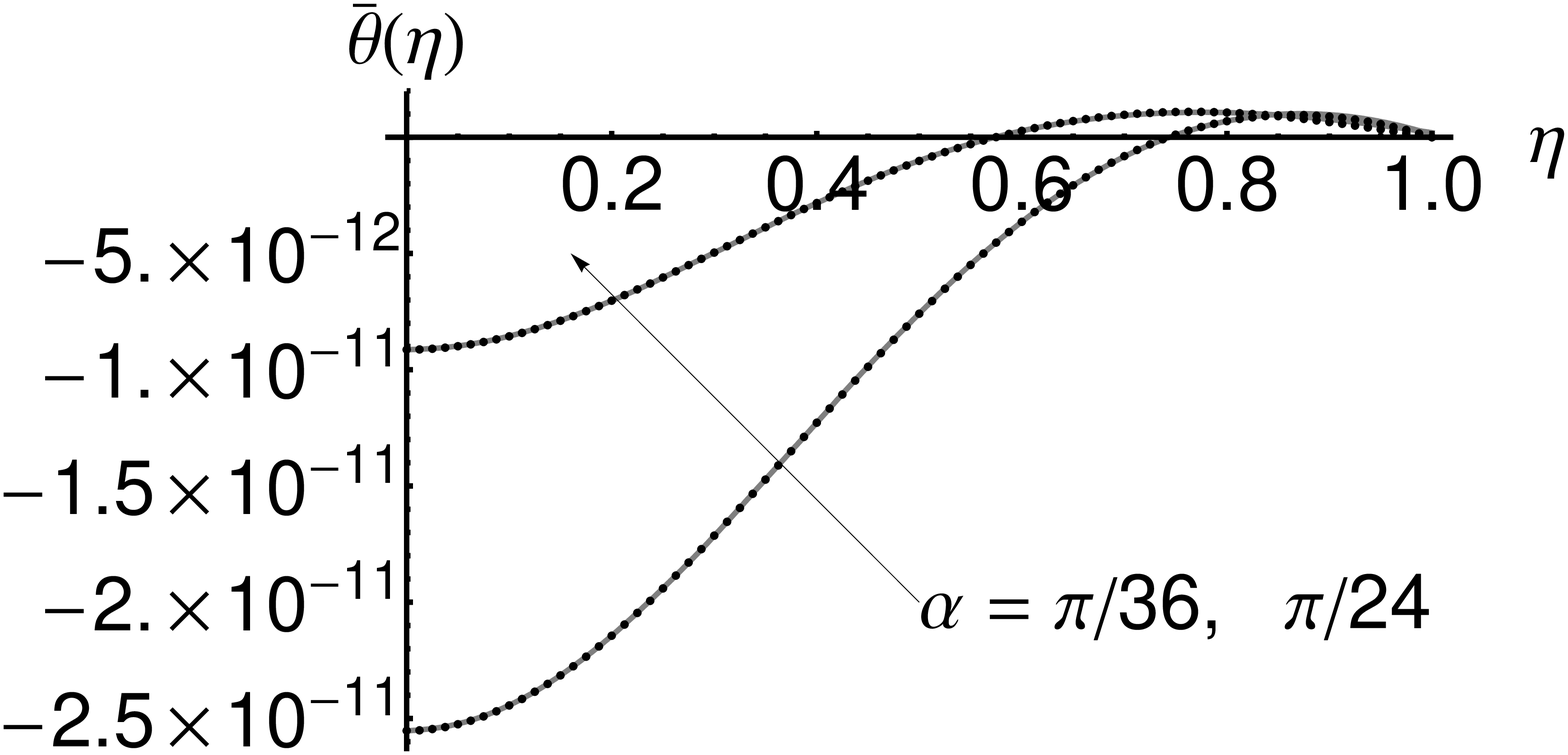}& \includegraphics[width=2.15in]{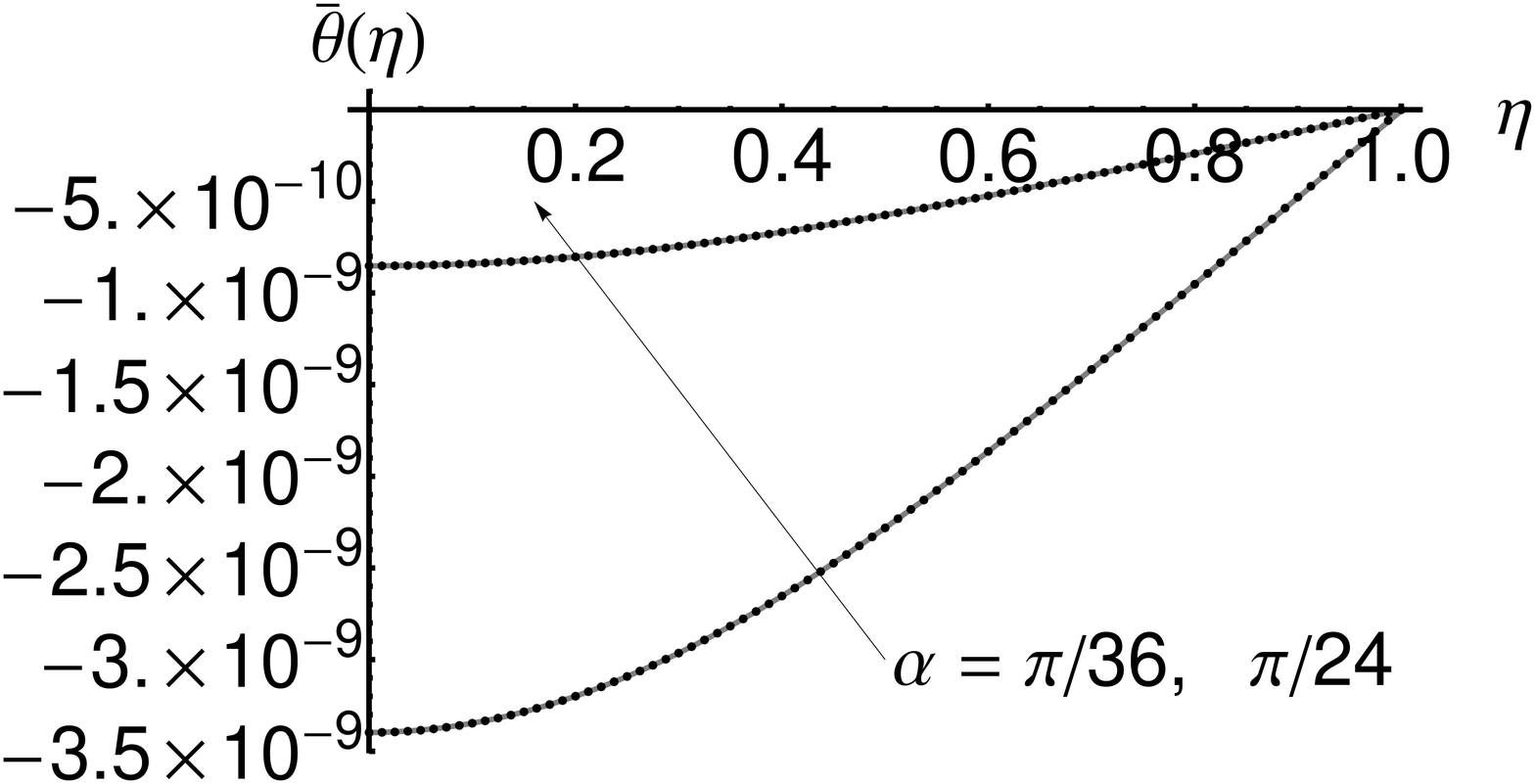}  \\
Fig. 10 Thermal profile  for $\alpha = \pi/36$,          &    Fig. 11 Thermal profile  for $\alpha = \pi/36$,     \\
and $\alpha = \pi/24$, $H=0$, $Re=50$,  &      and $\alpha = \pi/24$, $H=250$, $Re=50$, \\
Eqs. (\ref{JefferyHamel51}) and (\ref{JefferyHamel59}):   &   Eqs. (\ref{JefferyHamel53}) and (\ref{JefferyHamel61}): \\
------  numerical solution,      &    ------  numerical solution,       \\
........ OHPM solution &  
........ OHPM solution  \\
\end{tabular}\\
\end{tabular}\\
\begin{tabular}[!t]{c}
\begin{tabular}[!t]{c c}
  \includegraphics[width=2.15in]{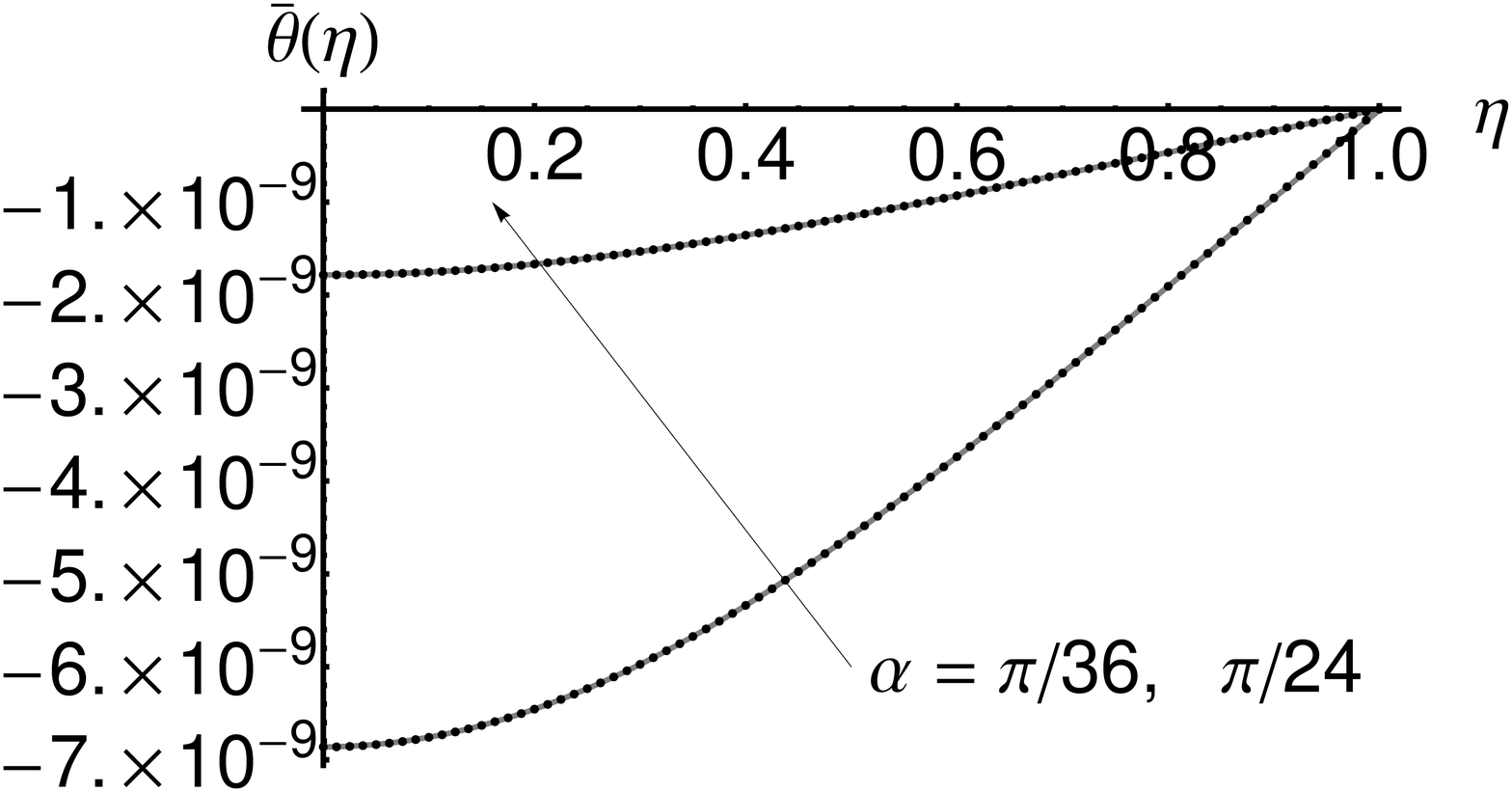}& \includegraphics[width=2.15in]{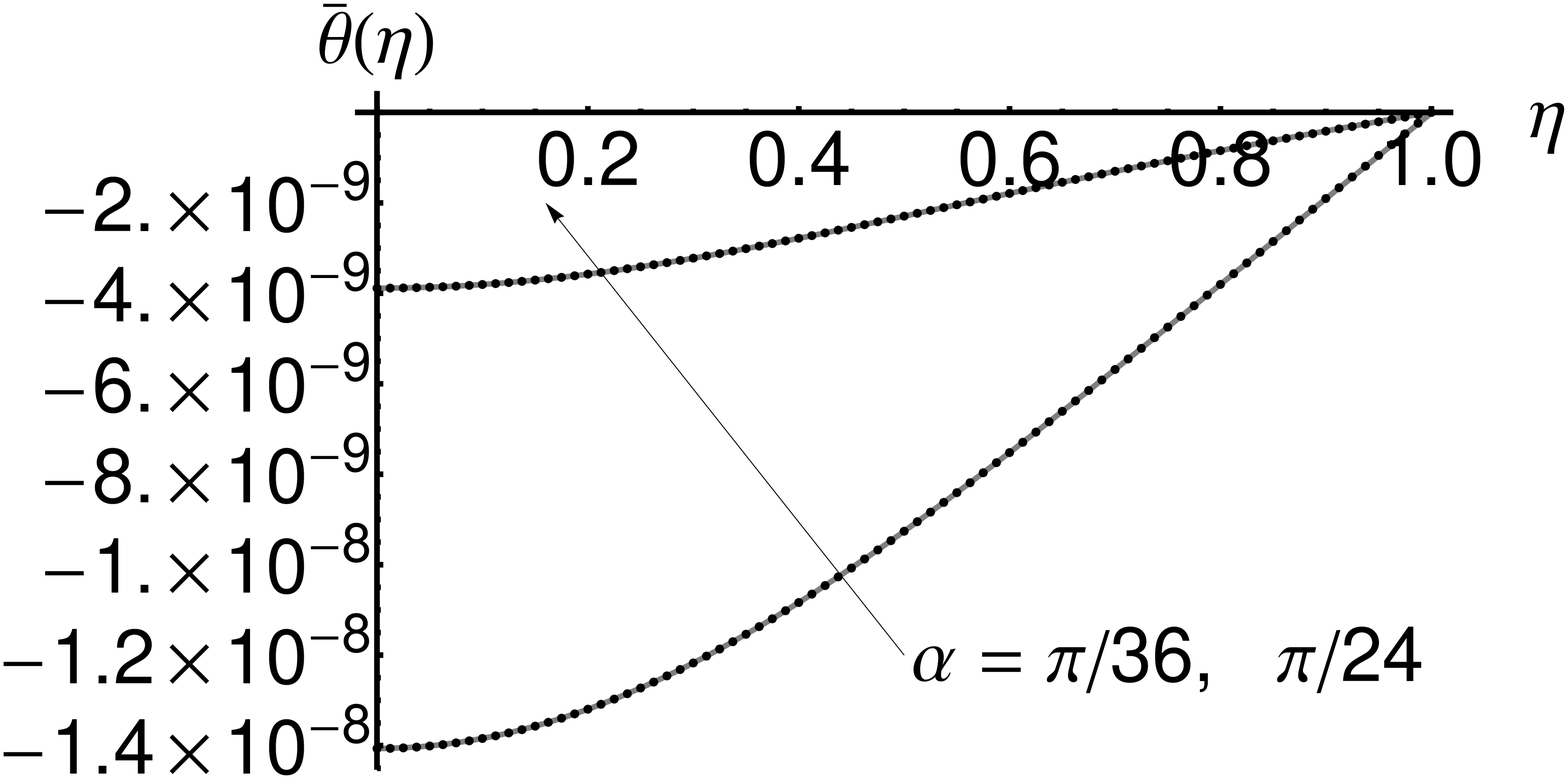}  \\
Fig. 12 Thermal profile  for $\alpha = \pi/36$,          &    Fig. 13 Thermal profile  for $\alpha = \pi/36$,     \\
and $\alpha = \pi/24$, $H=500$, $Re=50$,  &      and $\alpha = \pi/24$, $H=1000$, $Re=50$, \\
Eqs. (\ref{JefferyHamel55}) and (\ref{JefferyHamel63}):   &   Eqs. (\ref{JefferyHamel57}) and (\ref{JefferyHamel65}): \\
------  numerical solution,      &    ------  numerical solution,       \\
........ OHPM solution &  
........ OHPM solution  \\
\end{tabular}\\
\end{tabular}

\end{document}